\documentclass{amsart}
\usepackage[all]{xy}
\usepackage{multicol}
\usepackage{array}
\usepackage{graphicx}
\usepackage{amsmath}
\usepackage{amsthm}
\usepackage{amssymb}
\usepackage{amsfonts}
\usepackage{color}
\usepackage{enumerate}
\usepackage{mathrsfs}
\usepackage{pb-diagram}
\usepackage{lamsarrow,pb-lams}
\usepackage[all]{xy}

\def\CC{{\mathcal C}}
\def\D{{\Delta}}

\def\Q{{\mathbb Q}}

\def\dim{\operatorname{dim}}
\def\ann{\operatorname{ann}}

\def\m{{\mathfrak{m}}}

\def\L{{\Lambda}}

\def\w{{\omega}}
\def\W{{\Omega}}

\def\D{{\Delta}}
\def\E{{\mathcal{E}}}
\def\F{{\mathbb{F}}}
\def\FF{{\mathcal{F}}}
\def\G{{\mathcal{G}}}
\def\H{{\mathcal{H}}}

\def\I{{\mathcal{I}}}

\def\O{{\mathcal{O}}}
\def\PP{{\mathcal{P}}}

\def\to{{\rightarrow}}
\def\lto{{\longrightarrow}}

\def\proj{{\twoheadrightarrow}}

\def\Pic{{\textrm{Pic}}}

\def\Spec{{\textrm{Spec }}}

\def\dim{{\textrm{dim }}}
\def\codim{{\textrm{codim}}}
\def\Coh{{\textrm{Coh}}}

\def\Pic{{\textrm{Pic}}}

\def\adj{{\textrm{adj}}}

\def\Supp{{\textrm{Supp }}}

\def\hocolim{{\textrm{hocolim}_{\rightarrow}}}
\def\holim{{\textrm{holim}_{\leftarrow}}}

\newtheorem{defi}{Definition}[section]
\newtheorem{prop}[defi]{Proposition}
\newtheorem{corol}[defi]{Corollary}
\newtheorem{thm}[defi]{Theorem}

\newtheorem{lemma}[defi]{Lemma}

\newenvironment{pf}{\textit{Proof.}}{$\square$}
\newenvironment{pf.}{\textit{Proof }}{$\square$}
\newenvironment{skpf}{\textit{Sketch of proof.}}{$\square$}

\theoremstyle{definition}

\newtheorem{rmk}{Remark}[section]

\title[Irreducible Theta Divisors of PPAV's are Strongly F-regular]{Irreducible Theta Divisors of Principally Polarized Abelian Varieties are Strongly F-regular}
\author{Alan Marc Watson}
\address{Department of Mathematics, University of Utah, Salt Lake City, Utah 84112}
\email{watson@math.utah.edu}
\date{March 9, 2016}

\begin{document}

\maketitle

\begin{abstract}
We study the birational geometry of irregular varieties and the singularities of Theta divisors of PPAV's in positive characteristic by applying recent generic vanishing results of Hacon and Patakfalvi. In particular, we prove that irreducible Theta divisors of principally polarized abelian varieties are strongly F-regular, which extends an old result of Ein and Lazarsfeld to fields of positive characteristic. In order to prove this, we formulate a positive characteristic analogue of another result of Ein and Lazarsfeld, to the effect that the Albanese image of a smooth projective variety of maximal Albanese dimension with vanishing holomorphic Euler characteristic is fibered by abelian subvarieties.
\end{abstract}

\tableofcontents

\section{Introduction}
The purpose of this paper is to apply recent generic vanishing results in positive characteristic due to Hacon and Patakfalvi ~\cite{hp13} to the study of the birational geometry of irregular varieties and the singularities of Theta divisors of principally polarized abelian varieties. Over fields of characteristic zero, seminal work of Ein and Lazarsfeld ~\cite{el97} applied generic vanishing techniques over the complex numbers to settle a number of questions concerning the geometry of irregular varieties. One of their main results states that irreducible Theta divisors on principally polarized abelian varieties have mild singularities:

\begin{thm}(c.f. ~\cite[Theorem 1]{el97})
Let $A$ be an abelian variety and let $\Theta\subset A$ be a principal polarization (i.e. an ample divisor such that $h^0(A,\O_A(\Theta))=1$). If $\Theta$ is irreducible, then it is normal and has rational singularities. \label{el97-main-theorem}
\end{thm} \medskip

The conclusion of the theorem is captured by the adjoint ideal of $\Theta$: given any log resolution $\mu:A' \to A$ of the pair $(A,\Theta)$ and writing $\mu^{\ast}\Theta=\Theta'+F$ with $\Theta'$ smooth and $F$ $\mu$-exceptional, one may define $\adj(A,\Theta)=\mu_{\ast}\O_{A'}(K_{A'/A}-F)$. Standard arguments show that $\adj(A,\Theta)=\O_A$ is equivalent to $\Theta$ being normal and having rational singularities (see section 9.3.E in ~\cite{laz04}). Bearing this in mind, Ein and Lazarsfeld's argument breaks into the following steps: let $X\to \Theta$ be a resolution of singularities.

\begin{enumerate}[(i)]
\item If $\Theta$ is irreducible, then $X$ is of general type. This relies on a classical argument due to Ueno (see ~\cite{mor00}, Theorem 3.7), characterizing the Itaka fibration and the Kodaira dimension of an irreducible subvariety of an abelian variety.

\item If $X$ is of general type, then $\chi(\w_X)>0$. More concretely, if $X$ is a smooth projective variety of maximal Albanese dimension and $\chi(X,\omega_X)=0$, then the image of the Albanese map is fibred by tori. In particular, this shows that if $X$ is birational onto its image under the Albanese map, then $X$ is not of general type.

\item Generic vanishing theorems and Nadel vanishing yield $\adj(\Theta)=\O_A \Longleftrightarrow \chi(X,\omega_X)>0$.
\end{enumerate} \medskip

Therefore if $\Theta$ is irreducible, its adjoint ideal must be trivial, and by the characterization described above, it must be normal an have rational singularities. Work of Abramovich ~\cite{abr95} shows that the statement in (i) remains valid in positive characteristic once an appropriate notion of Kodaira dimension is defined for possibly singular varieties. In this paper we provide positive characteristic analogues of items (ii) and (iii). \medskip

The only known results in this direction are due to Hacon ~\cite{hac11}, where he proved that for principally polarized abelian varieties $(A,\Theta)$ over algebraically closed fields of positive characteristic, the pair $(A,\Theta)$ is a limit of strongly F-regular pairs. More precisely:

\begin{thm}(c.f. ~\cite[Theorem 1.1]{hac11})
Let $(A,\Theta)$ be a principally polarized abelian variety over an algebraically closed field of characteristic $p>0$. If $D\in |m\Theta|$, then $\left(A,\frac{1-\epsilon}{m}D\right)$ is strongly F-regular for any rational number $0<\epsilon<1$
\end{thm} \medskip

We summarize briefly our main results. The arguments employed in the proofs bear a strong resemblance to their characteristic zero analogues, albeit plenty of technicalities arise. Not only is resolution of singularities unavailable in general, but also generic vanishing for canonical sheaves is known to fail in positive characteristic (c.f. ~\cite{hk12}). Nevertheless, recent work of Hacon and Patakfalvi ~\cite{hp13} provides strong generic vanishing statements for objects arising from Cartier modules (see section 2.4 for precise statements): given a coherent Cartier module $\Omega_0\in \Coh(A)$, the traces of the Frobenius iterates yield an inverse system $$\cdots \to F_{\ast}^e\W_0 \to F_{\ast}^{e-1}\W_0 \to \cdots$$ and denoting by $\W=\varprojlim F_{\ast}^e\W_0$ its inverse limit, there exists a closed subset $Z\subset \hat{A}$ such that $H^i(A,\W\otimes P_{\alpha}^{\vee})=0$ for every $i>0$ and every $\alpha\in Z$ such that $p^e\alpha\notin Z$ for all $e>>0$ (c.f. Corollary 3.3.1 in ~\cite{hp13}). This grounds on the following Theorem, which is the main result in ~\cite{hp13}.

\begin{thm}(c.f. ~\cite[Theorem 3.1.1 and Lemma 3.1.2]{hp13})
Let $k$ be an algebraically closed field of characteristic $p>0$ and
$A$ be an abelian variety over $k$. Let $\{\Omega_e\}$ be Cartier module on $A$. If for any sufficiently
ample line bundle $L$ on $\hat{A}$ and any $e\gg 0$,
$H^i(A,\Omega_e\otimes \hat{L}^\vee)=0$ for all $i>0$, then the
complex\footnote{Here $S_{A,\hat{A}}$ denotes the Fourier-Mukai functor with kernel given by the Poincar{\'e} bundle of $A\times \hat{A}$} $$\Lambda=\hocolim RS_{A,\hat{A}}(D_A(\Omega_e))$$ is a
quasi-coherent sheaf in degree 0, i.e.,
$\Lambda=\mathcal{H}^0(\Lambda)$.
\end{thm} \medskip

This result is generalized further in ~\cite{wz14}, where a notion of M-regularity in positive characteristic is also introduced. Concretely, one has the following:

\begin{thm}[c.f. Theorem 4.2 in ~\cite{wz14}]
Let $A$ be an abelian variety and
$\{\Omega_e\}$ be a GV-inverse system of coherent sheaves on $A$
such that $\{\Omega_e\}$ is M-regular, in the sense that $\mathcal{H}^0(\Lambda)$ is torsion-free. Then for any scheme-theoretic point $P\in A$, if $\dim P\geqslant i$, then $P$ is not in the support of $$Im (R^i\hat{S}(\Omega) \to R^i\hat{S}(\Omega_e))$$ for any $e$.
\end{thm} \medskip

Our main technical result is a partial converse to the previous theorem: the presence of torsion in $\H^0(\L)$ induces the following non-vanishing statement:

\begin{thm}[c.f. Theorem 3.1]
Let $\{\Omega_e\}$ be an inverse system of coherent sheaves on a
$g$-dimensional abelian variety satisfying the Mittag-Leffler
condition and let $\Omega=\varprojlim \Omega_e$. Let
$\Lambda_e=RS_{A,\hat{A}}(D_A(\Omega_e))$ and $\Lambda=\hocolim\Lambda_e$.
Suppose that $\{\W_e\}$ is a GV-inverse system, in the sense that $\mathcal{H}^i(\Lambda)=0$ for any $i\neq 0$. If $\H^0(\L)$ is has a torsion point $P$ of dimension $g-k$, then the maps $$\varprojlim \left(R^kS_{A,\hat{A}}(\W_e)\otimes k(P)\right) \to R^kS_{A,\hat{A}}(\W_e)\otimes k(P)$$ are non-zero for every $e>>0$.
\end{thm}

Using this, we can derive a fibration statement similar to that of Ein and Lazarsfeld:

\begin{thm}[c.f. Theorem 4.2]
Let $X$ be a smooth projective variety of maximal Albanese dimension and denote by $a:X\to A$ the Albanese map. Let $g=\dim A$. Consider the inverse system $\{\W_e:=F_{\ast}^eS^0a_{\ast}\w_X\}_e$ and denote $\W=\varprojlim \W_e$. Define $\L_e=RS_{A,\hat{A}}D_A(\W_e)$ and assume that the sheaf $\H^0(\L)=\varinjlim \H^0(\L_e)$  has torsion. Then the image of the Albanese map is fibered by abelian subvarieties of $\hat{A}$.
\end{thm} \medskip

An identical argument to the one we employ to prove the previous theorem also yields the following result \label{Main-theorem} describing the singularities of Theta divisors in positive characteristic.

\begin{thm}
Let $A$ be an ordinary abelian variety over an algebraically closed field of positive characteristic and let $\Theta$ be an irreducible Theta divisor. Then $\Theta$ is strongly F-regular. \label{Main-theorem}
\end{thm} \medskip

Over fields of positive characteristic, work of Smith and Hara (c.f. \cite{smi97}, \cite{har98})  shows that F-rationality is the positive characteristic counterpart to rational singularities, and the former is implied by strong F-regularity, so in this sense Theorem \ref{Main-theorem} is stronger than one might expect. \medskip

This paper is structured as follows. We start by recording all the background results we need in section 2: in 2.1 we recall the main definitions and some useful properties of the Fourier-Mukai transform in the context of abelian varieties, in 2.2 we record the relevant definitions of F-singularities; in 2.3 we record results of Pink and Roessler characterizing subvarieties of abelian varieties in positive characteristic; in 2.4 we outline a few useful facts concerning inverse systems that will ease the exposition of the proofs and in 2.5 we collect the generic vanishing statements in positive characteristic that will be needed in the sequel. Sections 3 and 4 constitute the technical core of the paper: section 3 contains the proof of the non-vanishing statement in the presence of torsion of $\H^0(\L)$ and in section 4 we generalize Ein and Lazarsfeld's fibration statement. Finally in section 5 we present the proof of Theorem \ref{Main-theorem} on the singularities of Theta divisors. We start with the case of simple abelian varieties in section 5.1, which is a simple computation following easily from the results in ~\cite{hp13} that does not require the arguments from sections 3 and 4. The general case is more involved and is presented section 5.2.

\subsection*{Acknowledgements} The author would like to thank
his advisor Christopher Hacon for suggesting this problem, for invaluable discussions
and for sharing an early draft of \cite{hp13}.

\section{Preliminaries}

\subsection{Derived categories and Fourier-Mukai transforms}

Let $A$ be a g-dimensional abelian variety, denote by $\hat{A}=\Pic^0(A)$ its dual and let $\FF\in \Coh(A)$. Let $\PP\in \Pic(A\times \hat{A})$ be the Poincar{\'e} bundle and denote and consider the usual Fourier-Mukai functors: $$RS_{A,\hat{A}}^{\PP}:D(A) \to D(\hat{A}), \qquad RS_{A,\hat{A}}^{\PP}(\bullet)=Rp_{\hat{A}\ast}(p_A^{\ast}(\bullet)\otimes \PP)$$ even though we will most often omit $\PP$ from the notation and simply write $RS_{A,\hat{A}}(\bullet)$. \medskip

We start by stating Mukai's inversion theorem in the derived category of quasi-coherent sheaves:

\begin{thm}[~\cite{muk81}]
If $[-g]$ denotes the rightwise shift by $g$ places and $-1_A$ is the inverse on $A$, the following equalities hold on $D_{qc}(A)$ and $D_{qc}(\hat{A})$ $$RS_{\hat{A},A} \circ RS_{A,\hat{A}}= (-1_A)^{\ast}[-g], \qquad RS_{A,\hat{A}} \circ RS_{\hat{A},A}= (-1_{\hat{A}})^{\ast}[-g]$$
\end{thm} \medskip

We will also be using the following two results:

\begin{lemma}[~\cite{muk81}, Proposition 3.8]
The Fourier-Mukai transform commutes with the dualizing functor in $D_{qc}(\hat{A})$ up to inversions and shifts, namely $$D_A \circ RS_{\hat{A},A} \simeq \left((-1_{A})^{\ast} \circ RS_{\hat{A},A} \circ D_{\hat{A}}\right)[g]$$
\end{lemma} \medskip

\begin{lemma}(c.f, ~\cite[Lemma 3.4]{muk81})
Let $\phi:A\to B$ be an isogeny between abelian varieties and denote by $\hat{\phi}:\hat{B} \to \hat{A}$ the dual isogeny. Then the following equalities hold on $D_{qc}(B)$ and $D_{qc}(A)$ respectively. $$\phi^{\ast} \circ RS_{\hat{B},B} \simeq RS_{\hat{A},A} \circ \hat{\phi}_{\ast}, \qquad \phi_{\ast} \circ RS_{\hat{A},A} \simeq RS_{\hat{B},B} \circ \hat{\phi}^{\ast}$$ In particular, this holds for the (e-th iterate) Frobenius map $F^e$ and its dual isogeny, namely the Verschiebung map $V^e=\hat{F}^e$. \label{Mukai-vs-isogenies}
\end{lemma}

We will also be using the following simple remark.

\begin{lemma}(c.f. ~\cite[Exercise 5.12]{huy06})
Let $\pi:B\to A$ be a morphism between abelian varieties and let $\PP$ be a locally free sheaf on $A\times \hat{A}$. Denote $\PP_{\pi}=(\pi \times 1_{\hat{A}})^{\ast}(\PP)$. Then $$S_{\PP_{\pi}} \simeq S_{\PP} \circ \pi_{\ast}$$ \label{Mukai-vs-push-forward}
\end{lemma} \medskip

We next record the notions of homotopy limits and colimits in the derived category. Given a direct system of objects $\CC_i\in D(A)$ $$\CC_1 \to \CC_2 \to \ldots$$ its homotopy colimit $\hocolim \CC_i$ is defined by the triangle $$\oplus \CC_i \lto \oplus \CC_i \lto \hocolim \CC_i \stackrel{[+1]}{\lto}$$ where the first map is the homomorphism given by $id-shift$ where $shift:\oplus \CC_i \to \oplus \CC_i$ is given on $\CC_i$ by the composition $\CC_i \to \CC_{i+1} \hookrightarrow \oplus \CC_j$. \medskip

Given an inverse system of objects $\CC_i \in D_{qc}(X)$ $$\CC_1 \longleftarrow \CC_2 \longleftarrow \cdots$$ its homotopy limit $\holim \CC_i$ is given by the triangle $$\holim \CC_i \lto \prod \CC_i \lto \prod \CC_i \stackrel{+1}{\lto}$$ where the map between products is $\prod(id-shift)$ and where by product we mean product of chain complexes as opposed to the product inside $D_{qc}(X)$. \medskip

Note that if $\CC_i$ are coherent sheaves, then $\hocolim \CC_i = \varinjlim \CC_i$. \medskip

If $X$ is an n-dimensional variety over a field $k$ and $\w_X^{\bullet}$ denotes its dualizing complex, so that $\H^{-\dim X}(\w_X^{\bullet}) \simeq \w_X$, we define the dualizing functor $D_X$ on $D_{qc}(X)$ as $D_X(\FF)= R\H om(\FF,\w_X^{\bullet})$. In this context, Grothendieck duality reads as follows:

\begin{thm}
Let $f:X\to Y$ be a proper morphism of quasi-projective varieties over a field $k$. Then for any complex $\FF\in D_{qc}(X)$ we have an isomorphism $$Rf_{\ast}D_X(\FF) \simeq D_Y Rf_{\ast}(\FF)$$ Assuming that $X$ and $Y$ are smooth, then we equivalently have that for any $\FF \in D_{qc}(X)$ and $\E\in D_{qc}(Y)$, if $\w_f=\w_X\otimes f^{\ast}\w_Y$ denotes the relative dualizing sheaf, there is a functorial isomorphism $$Rf_{\ast} R\H om(\FF,Lf^{\ast}(\E) \otimes \w_f[\dim X-\dim Y]) \simeq R\H om (Rf_{\ast}\FF,\E)$$ \label{grothendieck-duality}
\end{thm}

\subsection{F-singularities and linear subvarieties of abelian subvarieties}

In this section we recall the basic notions from the theory of F-singularities following \cite{sch12} and ~\cite{bst12}. Let $X$ be a separated scheme of finite type over an F-finite perfect field of characteristic $p>0$. A variety is a connected reduced equidimensional scheme over $k$. We denote the canonical sheaf of $X$ by $\w_X=\H^{-\dim X}(\w_X^{\bullet})$, where $\w_X^{\bullet}=\eta^{\ast}k$ is the dualizing complex of $X$ and $\eta:X \to k$ is the structural map. If $X$ is normal, a \textit{canonical divisor} on $X$ is any divisor $K_X$ such that $\w_X \simeq \O_X(K_X)$. \medskip

By a pair $(X,\D)$ we mean the combined information of a normal integral scheme $X$ and an effective $\Q$-divisor $\D$. Denote by $F^e:X \to X$ the e-th iterated absolute Frobenius, where the source has structure map $\eta\circ F^e:X\to k$. Since $(F^e)^!\w_X^{\bullet} = (F^e)^!\eta^!k=\eta^! (F^e)^!k=\eta^!k=\w_X^{\bullet}$. In general for a finite morphism $f:X\to Y$, a coherent sheaf $\FF$ on $X$ and a quasi-coherent sheaf $\G$ on $Y$, we have the duality $\H om(f_{\ast}\FF,\G) \simeq f_{\ast}\H om(\FF, f^!\G)$, so the identity $\w_X^{\bullet} \to \w_X^{\bullet} \simeq (F^e)^!\w_X^{\bullet}$ yields a trace map $F_{\ast}^e\w_X^{\bullet} \to \w_X^{\bullet}$ and taking cohomology we obtain $\Phi^e:F_{\ast}^e\w_X \to \w_X$. Given a variety $X$, the \textit{parameter test submodule} $\tau(\w_X)$ of $X$ is the unique smallest $\O_X$-submodule $M\subseteq \w_X$, non-zero on any component of $X$, such that $\Phi^1(F_{\ast}M) \subseteq M$. \medskip

Assume that $(X,\Delta)$ is a pair such that $K_X+\D$ is $\Q$-Cartier with index not divisible by $p$. Choose $e>0$ such that $(p^e-1)(K_X+\D)$ is Cartier and define the line bundle $\mathcal{L}_{e,\D}=\O_X((1-p^e)(K_X+\D))$. By ~\cite{sch09}, there is a canonically determined map $\phi_{e,\Delta}:F_{\ast}^e\mathcal{L}_{e,\D} \to \O_X$. We define the \textit{test ideal} $\tau(X,\D)$ of the pair $(X,\D)$ to be the smallest non-zero ideal $J\subseteq \O_X$ such that $$\phi_{e,\D}(F_{\ast}^e(J\cdot \mathcal{L}_{e,\D})) \subseteq J.$$ Similarly one defines the \textit{non-F-pure ideal} $\sigma(X,\D)$ of $(X,\D)$ to be the the largest such ideal $J\subseteq \O_X$. \medskip

Ever since Hochster and Huneke introduced test ideals and tight closure theory in ~\cite{hh90}, deep connections have been established between the classes of singularities defined in terms of Frobenius splittings and those arising within the minimal model program. For instance, a normal domain $(R,\m)$ of characteristic $p>0$ is said to be F-pure if the inclusion induced by the Frobenius $R\hookrightarrow F_{\ast}^eR\equiv R^{1/p^e}$ splits for every $e$. Similarly, a pair $(R,\D)$ is said to be F-pure if the inclusion $R\hookrightarrow R^{1/p^e} \hookrightarrow R\left(\lceil(p^e-1)\D\rceil\right)^{1/p^e}$ splits for every $e$ and it was shown in ~\cite{hw02} that F-pure pairs are the analogues of log canonical pairs in characteristic zero, in the sense that if $(X,\D)$ is a log canonical pair, then its reduction mod $p$ $(X_p,\D_p)$ is F-pure for all $p>>0$. \medskip

In this paper we shall be concerned with the two classes of F-singularities that we define next. We will be recording the original definition in terms of Frobenius splittings and we will then state their description in terms of test ideals that will be used in the sequel.

\begin{defi}
\begin{enumerate}[(i)]
\item A pair $(X=\Spec R,\D)$ is \textit{strongly F-regular} if for every non-zero element $c\in R$, there exists $e$ such that the map $R\hookrightarrow R^{1/p^e} \hookrightarrow R((p^e-1)\D)^{1/p^e}$ that sends $1\mapsto c^{1/p^e} \mapsto c^{1/p^e}$ splits as an $R$-module homomorphism.
\item A reduced connected variety $X$ is F-rational if it is Cohen-Macaulay and there is no non-zero submodule $M\subsetneq \w_R$ such that the Grothendieck trace map $\Phi_X:F_{\ast}^e\w_X \to \w_X$ satisfies $\Phi(F_{\ast}^e M)\subseteq M$.
\end{enumerate}
\end{defi} \medskip

Strongly F-regular pairs are the analog of log terminal pairs in characteristic zero (c.f. ~\cite{hw02}) and F-rational varieties are the analogue of varieties with rational singularities (c.f. ~\cite{smi97}). The notion of strong F-regularity is also captured by the test ideal, as the following well-known result shows.

\begin{lemma}(c.f. ~\cite[Proposition 2.4]{hw02})
A pair $(X,\D)$ is strongly F-regular if, and only if, $\tau(X,\D)=\O_X$.
\end{lemma}

Assume that $(X,\D)$ is a pair, where $X$ is a normal proper variety over an algebraically closed field of characteristic $p>0$ and $\D\geq 0$ is a $\Q$-divisor such that $K_X+\D$ is $\Q$-Cartier with index not divisible by $p$. The map $\phi_{\Delta}^e:F_{\ast}^e\mathcal{L}_{e,\D} \to \O_X$ defined in ~\cite{sch09} restricts to surjective maps $$F_{\ast}^e(\sigma(X,\D) \otimes \mathcal{L}_{e,\D}) \lto \sigma(X,D), \qquad F_{\ast}^e(\tau(X,\D) \otimes \mathcal{L}_{e,\D}) \lto \tau(X,D).$$

The power of vanishing theorems in characteristic zero relies on the fact that they allow us to lift global sections of adjoint bundles. The full space of global sections is not so well behaved in positive characteristic, so one instead focuses on a subspace of it that is stable under the Frobenius action. \medskip

If $M$ is any Cartier divisor, one thus defines the subspace $S^0(X, \tau(X,\D) \otimes \O_X(M))$ as

\begin{multline*} S^0(X, \tau(X,\D) \otimes \O_X(M)) \\ := \bigcap_{n\geq0} Im\left( H^0(X,F_{\ast}^{ne}\tau(X,\D) \otimes \mathcal{L}_{ne,\D}(p^{ne}M)) \lto H^0(X,\tau(X,\D) \otimes \O_X(M)) \right) \\ \subseteq H^0(X,\O_X(M)) \end{multline*}

Among the many applications of these subspaces, for instance, they can be used to prove global generation statements: concretely, suppose that $X$ is a $d$-dimensional variety of characteristic $p>0$ and that $\D$ is a $\Q$-divisor such that $K_X+\D$ is $\Q$-Cartier with index not divisible by $p$. It was shown in ~\cite{sch09} that if $L$ and $M$ are Cartier divisors such that $L-K_X-\D$ is ample and $M$ is ample and globally generated, then the sheaf $\tau(X,\D)\otimes \O_X(L+nM)$ is globally generated for all $n\geq d$ by $S^0(X, \tau(X,\D) \otimes \O_X(L+nM))$. \medskip

\subsection{The Frobenius morphism on Abelian varieties}

Throughout this paper, $A$ will denote an abelian variety of dimension $g$ over a field $k$ and $\hat{A}=\Pic^0(A)$ will denote the dual abelian variety

\begin{lemma}(c.f. ~\cite[Proposition 2.13]{hp13})
For a g-dimensional abelian variety $A$ over a field $k$, the following conditions are equivalent.
\begin{enumerate}[(i)]
\item There are $p^g$ p-torsion points.
\item The Frobenius action $H^g(A,\O_A) \to H^g(A,\O_A)$ is bijective
\item The Frobenius action $H^i(A,\O_A) \to H^i(A,\O_A)$ is bijective for all $0\leq i \leq g$
\item $S^0(A,\w_A) = H^0(A, \w_A)$
\end{enumerate}
\end{lemma} \medskip

If any of these equivalent conditions is satisfied we say that $A$ is \textit{ordinary}. Given an isogeny $\varphi:A \to B$ between abelian varieties of dimension $g$, $A$ is ordinary if and only if $B$ is ordinary. Given a surjective morphism $\varphi:A \to B$ of abelian varieties, if $A$ is ordinary then so is $B$ (see Lemmas 2.14 and 2.14 in ~\cite{hp13}). \medskip

We finally record a characterization of linear subvarieties of abelian varieties following ~\cite{pr03}. Let $A$ be an abelian variety endowed with an isogeny $\varphi:A \to A$. We say that $A$ is pure of positive weight if there exist integers $r,s>0$ such that $\varphi^s=F_{p^r}$ for some model of $A$ over $\F_{p^r}$. If $A$ is defined over a finite field, we say $A$ is \textit{supersingular} if and only if it is pure of positive weight for the isogeny given by multiplication by $p$; in general, we say that $A$ is supersingular if it is isogenous to a supersingular variety defined over a finite field. We say that $A$ \textit{has no supersingular factors} is there exist no non-trivial homomorphism  to an abelian variety which is pure of positive weight for the isogeny given by multiplication by $p$. One sees that $A$ has no supersingular factors if there does not exist a non-trivial homomorphism to a supersingular abelian variety. In particular, if $A$ is an ordinary abelian variety, it follows from the observations in the previous paragraph that $A$ has no supersingular factors (see Lemma 2.16 in ~\cite{hp13}). \medskip

The following result of Pink and Roessler characterizing linear subvarieties of abelian varieties will be crucial in our proof:

\begin{thm}(c.f. ~\cite[Theorem 2.2]{pr03})
Let $A$ be an abelian variety over a field $K$ of characteristic $p>0$ and let $X\subset A$ be a reduced closed subscheme $p(X) \subset X$, where $p$ denotes the isogeny given by multiplication by $p$. If $A$ has no supersingular factors, then all the maximal dimensional irreducible components of $X$ are completely linear (namely, torsion translates of subabelian varieties).
\end{thm} \medskip

\subsection{Generalities on inverse systems and spectral sequences}

\textbf{Mittag-Leffler inverse systems}. We start by recording a few results that will be used in the sequel, most of which are taken directly from ~\cite{har78}. Recall that a sheaf is countably quasi-coherent if it is quasi-coherent and locally countably generated. Also recall that an inverse system of coherent sheaves $\{\Omega_e\}$ is said to satisfy the Mittag-Leffler condition if for any $e\geq 0$, the image of $\Omega_{e'}\to \Omega_{e}$ stabilizes for $e'$ sufficiently large. The inverse limit functor is always left exact in the sense that if $$\xymatrix{0 \ar[r] & \mathcal{F}_e \ar[r] \ar[d] & \mathcal{G}_e \ar[r] \ar[d] & \mathcal{H}_e \ar[r] \ar[d] & 0 \\ 0 \ar[r] & \mathcal{F}_{e-1} \ar[r] & \mathcal{G}_{e-1} \ar[r] & \mathcal{H}_{e-1} \ar[r] & 0}$$ is an exact sequence of inverse systems, then $$0 \to \varprojlim\mathcal{F}_e \to \varprojlim\mathcal{G}_e \to \varprojlim\mathcal{H}_e$$ is exact in the category of quasi-coherent sheaves. By a theorem of Roos (c.f. Proposition I.4.1 in ~\cite{har78}), the right derived functors $R^i\varprojlim$ are $0$ for $i\geqslant 2$. Hence, we have a long exact sequence $$0 \to \varprojlim\mathcal{F}_e \to \varprojlim\mathcal{G}_e \to \varprojlim\mathcal{H}_e \to R^1\varprojlim\mathcal{F}_e \to R^1\varprojlim\mathcal{G}_e \to R^1\varprojlim\mathcal{H}_e \to 0.$$

We start by recording a characterization of the Mittag-Leffler condition in terms of the first right-derived inverse limit.

\begin{lemma}(c.f. ~\cite[Proposition I.4.9]{har78})
Let $\{\W_e\}_e$ be an inverse system of countably quasi-coherent sheaves on a scheme $X$ of finite type. Then the following conditions are equivalent:
\begin{enumerate}[(i)]
\item $\{\W_e\}_e$ is satisfies the Mittag-Leffler condition.
\item $R^1\varprojlim \W_e=0$
\item $R^1\varprojlim \W_e$ is countably quasi-coherent.
\end{enumerate}
\label{ML-charact}
\end{lemma}\medskip

The following is basic result about the cohomology of an inverse system of sheaves:

\begin{prop}(c.f. ~\cite[Theorem I.4.5]{har78})
Let $\{\Omega_e\}$ be an inverse system of coherent sheaves on a variety $X$. Let $T$ be a functor on $D(X)$ which commutes with arbitrary direct products. Suppose that $\{\Omega_e\}$ satisfies the Mittag-Leffler condition. Then for each $i$, there is an exact sequence $$0\to R^1\varprojlim R^{i-1}T(\Omega_e)\to R^iT(\varprojlim \Omega_e)\to \varprojlim R^iT(\Omega_e) \to 0.$$ In particular, if for some $i$, $\{R^{i-1}T(\Omega_e)\}$ satisfies the Mittag-Leffler condition, then $R^iT(\varprojlim \Omega_e)\cong \varprojlim R^iT(\Omega_e)$ (by Lemma \ref{ML-charact}). \label{inverse-limits-commute-functor}
\end{prop} \medskip

We will be applying this theorem to the push-forward $f_{\ast}$ under a proper morphism of schemes. We finally record a standard statement about the commutation of inverse limits and tensor products. Recall that a sheaf is countably quasi-coherent if it is quasi-coherent and locally countably generated. Then one has the following:

\begin{lemma}(see ~\cite[Proposition 4.10]{har78})
Let $\{\FF_e\}_e$ be an inverse system of countably quasi-coherent sheaves on a scheme $X$ of finite type and let $E$ be a flat $\O_X$-module. Consider the natural map $$\alpha:(\varprojlim \FF_e)\otimes E \rightarrow \varprojlim (\FF_e \otimes E)$$ If $\varprojlim \FF_e$ is countably quasi-coherent then $\alpha$ is injective and if furthermore $R^1\varprojlim \W_e$ is countably quasi-coherent, then $\alpha$ is surjective. In particular, if $\{\FF_e\}$ is an inverse system of coherent sheaves satisfying the Mittag-Leffler condition on a scheme $X$  with generic point $w$, then there is an isomorphism $$\left(\varprojlim_e \FF_e \right) \otimes k(w) \lto \left(\varprojlim_e \FF_e \otimes k(w) \right)$$
\label{inverse-limit-tensor-product-commute}
\end{lemma} \medskip

\textbf{Inverse systems of convergent spectral sequences}. The following observation is taken from ~\cite{car08}. Let $\{E(n)\}$ be an inverse system of spectral sequences with morphisms of spectral sequences $E(n) \to E(n-1)$ and consider the tri-graded abelian groups $E_{p,q}^r=\varprojlim_n E_{p,q}^r(n)$, with differentials given by the inverse limits of the differentials in the $E(n)$. Concretely, if $d^r(n)$ is the r-th differential in $E(n)$ and $x(n)\in E_{p,q}^r$, then the r-th differential $d^r$ in the limit sequence $E_{p,q}^r$ is given by $d^r(x(n))=d^r(n)(x(n))$. The resulting object is a spectral sequence provided that $H(E_{p,q}^r,d^r)=E_{p,q}^{r+1}$, which is in turn equivalent to showing that $$H(\varprojlim_n E_{p,q}^r(n),d^r)=\varprojlim_n H(E_{p,q}^r(n),d^r)$$ and this is precisely the statement of Proposition \ref{inverse-limits-commute-functor} above (for the functor of global sections). \medskip

Note that, in particular, if the terms $E_{p,q}^r$ are all finite-dimensional vector spaces, the hypotheses of Proposition \ref{inverse-limits-commute-functor} hold and all the $\varprojlim^{(1)}$ terms are 0. Besides, given that the inverse limit of the spectral sequences is again a spectral sequence and provided that every spectral sequence in the inverse system is bounded and convergent, one observes that the limit spectral sequence is also convergent: for fixed $p,q$, there is a fixed $N$ such that $E_{p,q}^{\infty}(n)=E_{p,q}^N(n)$.\medskip

\textbf{Morphisms between spectral sequences}. We recall the definition of a spectral sequence from EGA III [$0_{III}$, 1.1]. Let $\mathscr{C}$ be an abelian category. A \textbf{(biregular) spectral sequence} $E$ on $\mathscr{C}$ consists of the following ingredients:
\begin{enumerate}
\item A family of objects $\{E^{p,q}_{r}\}$ in $\mathscr{C}$,
where $p,q,r\in \mathbb{Z}$ and $r\geqslant 2$, such that for any
fixed pair $(p,q)$, $E^{p,q}_r$ stabilizes when $r$ is sufficiently
large. We denote the stable objects by $E^{p,q}_\infty$.
\item A family of morphisms $d^{p,q}_r:E^{p,q}_r\to
E^{p+r,q-r+1}_r$ satisfying $$d^{p+r,q-r+1}_r\circ d^{p,q}_r=0.$$
\item A family of isomorphisms
$\alpha^{p,q}_r:\ker(d^{p,q}_r)/Im(d^{p-r,q+r-1}_r)\stackrel{\rightarrow}{\sim}
E^{p,q}_{r+1}$.
\item A family of objects $\{E^n\}$ in $\mathscr{C}$. For every
$E_n$, there is a bounded decreasing filtration $\{F^pE^n\}$ in the
sense that there is some $p$ such that $F^pE^n=E^n$ and there is
some $p$ such that $F^pE^n=0$.
\item A family of isomorphisms
$\beta^{p,q}:E^{p,q}_\infty\stackrel{\rightarrow}{\sim}F^pE^{p+q}/F^{p+1}E^{p+q}$.
\end{enumerate}
We say that the spectral sequence $\{E^{p,q}_r\}$ converges to $\{E^n\}$
and write $$E^{p,q}_2\Rightarrow E^{p+q}.$$

A morphism $\phi:E\to H$ between two spectral sequences on
$\mathscr{C}$ is a family of morphisms $\phi^{p,q}_r:E^{p,q}_r\to
H^{p,q}_r$ and $\phi^n:E^n\to H^n$ such that $\phi$ is compatible
with $d$, $\alpha$, the filtration and $\beta$. The following result is useful in order to obtain information about the limiting map $\phi^n$ from the maps $\phi_2^{p,q}$.

\begin{lemma}[see Lemma 2.15 in ~\cite{wz14}]
Let $$\xymatrix{ E_2^{i,j} \ar@2{->}[r] \ar^{\phi_2^{i,j}}[d] & E^{i+j} \ar^{\phi^{i+j}}[d] \\
H_2^{i,j} \ar@2{->}[r] & H^{i+j}}$$ be two spectral sequences with
commutative maps. Let $l$ and $a$ be integers. Suppose that
$E_2^{i,l-i}=0$ for $i<a$, $H_2^{i,l-i}=0$ for $i>a$ and
$\phi_2^{a,l-a}=0$. Then $\phi^l=0$. \label{spec-seq-zero-map}
\end{lemma} \medskip

\subsection{Generic vanishing in positive characteristic}

A smooth projective variety $X$ over an algebraically closed field is said to have maximal Albanese dimension if it admits a generically finite morphism to an abelian variety $X\to A$. Over fields of characteristic zero, the main tool that is employed when studying properties of varieties of maximal Albanese dimension is the generic vanishing theorem of Green and Lazarsfeld (~\cite{gl87}, ~\cite{gl91}). Even though it is shown in ~\cite{hk12} that the obvious generalization of this result to fields of positive characteristic if false, recent work of Hacon and Patakfalvi ~\cite{hp13} provides a weaker generic vanishing statement in positive characteristic which albeit necessarily weaker, is strong enough prove positive characteristic versions of Kawamata's celebrated characterization of abelian varieties. In this subsection we collect the results of ~\cite{hp13} that we shall be using throughout. \medskip

The following is the main theorem in ~\cite{hp13}:

\begin{thm}( c.f. ~\cite[Theorem 3.1, Lemma 3.2]{hp13})
Let $A$ be an abelian variety defined over an algebraically closed field of positive characteristic and let $\Omega_{e+1} \to \W_e$ be an inverse system of coherent sheaves on $A$.

\begin{enumerate}[(i)]
\item If for any sufficiently ample line bundle $L\in \Pic(\hat{A})$ and for any $e>>0$ we have $H^i(A,\W_e \otimes RS_{\hat{A},A}(L)^{\vee})=0$ for every $i>0$, then the complex $\Lambda=\hocolim RS_{A,\hat{A}}(D_A\W_e)$ (which in general is concentrated in degrees $[-g,\ldots,0]$), is actually a quasi-coherent sheaf concentrated in degree 0, namely $\Lambda=\H^0(\Lambda)$. Besides $\W=\varprojlim \W_e = \left((-1_A)^{\ast}D_A RS_{\hat{A},A}(\Lambda)\right)[g]$.

\item The condition in (i) is satisfied for coherent Cartier modules: if $F_{\ast}\W_0 \to \W_0$ is a coherent Cartier module and we denote $\W_e=F_{\ast}^e\W_0$, then for any ample line bundle $L\in \Pic(\hat{A})$ we have $$H^i(A,\W_e \otimes RS_{\hat{A},A}(L)^{\vee} \otimes P_{\alpha})=0, \quad \forall e>>0, \quad \forall i>0, \quad \forall \alpha\in \hat{A}$$
\end{enumerate} \label{generic-vanishing-char-p}
\end{thm} \medskip

From the above result and the cohomology and base change theorem one derives the following corollary:

\begin{corol}(c.f. ~\cite[Corollaries 3.5 and 3.6]{hp13})
With the same notations as above we have the following:

\begin{enumerate}[(i)]
\item For every $\alpha\in \hat{A}$ we have $\Lambda \otimes k(\alpha) \simeq \varinjlim H^0(A,\W_e \otimes \PP_{\alpha}^{\vee})^{\vee}$, and for every integer $e\geq0$, $\H^0(\Lambda_e) \otimes k(\alpha) \simeq H^0(A,\W_e \otimes \PP_{\alpha}^{\vee})^{\vee}$.

\item There exists a proper closed subset $Z\subset \hat{A}$ such that if $i>0$ and $p^ey\notin Z$ for all $e>>0$, then $\varinjlim H^i(A,\W_e \otimes \PP_{\alpha}^{\vee})^{\vee}=0$. Furthermore, if $W^i=\{\alpha\in \hat{A}, \quad \varprojlim H^i(A,\W_e \otimes \PP_{\alpha}^{\vee})^{\vee} \neq 0\}$, then $$W^i \subset Z'=\overline{\bigcup_{e\geq 0} \left([p_{\hat{A}}^e]^{-1}(Z)\right)_{red}}$$ where $pZ' \subset Z'$ If besides $\hat{A}$ has no supersingular factors, then the top dimensional components of $Z'$ are a finite union of torsion translates of subtori of $A$.
\end{enumerate} \label{GV-corollary1}
\end{corol} \medskip

We quote two more results from ~\cite{hp13} that will provide a simple proof of a special case of our main theorem:

\begin{prop}(c.f. ~\cite[Proposition 3.17]{hp13})
Let $A$ be an ordinary abelian variety and consider the same notations as above. Then each maximal dimensional irreducible component of the set $Z$ of points $\alpha \in \hat{A}$ such that the image of the natural map $$\H^0(\Lambda_0) \otimes \O_{\hat{A},\alpha} \lto \H^0(\Lambda) \otimes \O_{\hat{A},\alpha} \simeq \Lambda \otimes \O_{\hat{A},\alpha}$$ is non-zero, is a torsion translate of an abelian subvariety of $\hat{A}$ and $\Lambda \otimes \O_{\hat{A},\alpha} \neq 0$ if and only if $\PP_{\alpha}^e \in Z$. \label{GV-corollary2}
\end{prop} \medskip

\begin{prop}(c.f. ~\cite[Lemma 3.9, Corollary 3.10]{hp13})
Let $\W_0$ be a coherent sheaf on an abelian variety $A$ and assume that $F_{\ast}\W_0 \to \W_0$ is surjective. Then $\Supp \W = \Supp \W_0$, so that $\Supp \W$ is a closed subvariety. Let $\hat{B} \subset \hat{A}$ be an abelian subvariety such that $$V^0(\W_0)=\{\alpha\in \hat{A}: h^0(\W_0 \otimes \PP_{\alpha}) \neq 0\}$$ is contained in finitely many translates of $\hat{B}$. Then $t_x^{\ast}\W \simeq \W$ for every $x\in \widehat{\hat{A}/\hat{B}}$. In particular, $\Supp \W$ is fibered by the projection $A\to B$, namely $\Supp \W$ is a union of fibers of $A\to B$. \label{GV-corollary3}
\end{prop} \medskip

Note that, in particular, Proposition \ref{GV-corollary3} applies to the subvariety $$Z=\left\{\alpha\in \hat{A}: \quad Im \left(\H^0(\Lambda_0)\otimes \O_{\hat{A},\alpha} \lto \Lambda \otimes \O_{\hat{A},\alpha} \right) \neq 0\right\}$$ from Proposition \ref{GV-corollary2}. \medskip

Finally, grounding on the results in ~\cite{hp13}, part of Pareschi and Popa's generic vanishing theory was extended to positive characteristic in ~\cite{wz14}. The main result in that paper is the following:

\begin{thm} Let $A$ be an abelian variety. Let $\{\Omega_e\}$ be an inverse system of coherent sheaves on $A$ satisfying the Mittag-Leffler condition and let $\Omega=\varprojlim \Omega_e$. Let $\Lambda_e=R\hat{S}(D_A(\Omega_e))$ and $\Lambda=\hocolim\Lambda_e$. The following are equivalent:

\begin{enumerate}
\item[(1)] For any ample line bundle $L$ on $\hat{A}$, $H^i(A,\Omega\otimes \hat{L}^\vee)=0$ for any $i>0$.
\item[(1')] For any fixed positive integer $e$ and any $i>0$, the homomorphism $$H^i(A, \Omega\otimes \hat{L}^\vee) \to H^i(A,\Omega_e\otimes \hat{L}^\vee)$$ is 0 for any sufficiently ample line bundle $L$.
\item[(2)] $\mathcal{H}^i(\Lambda)=0$ for any $i\neq 0$.
\end{enumerate}

These conditions imply the following:

\begin{enumerate}
\item[(3)] For any scheme-theoretical point $P\in A$, if $\dim P>g-i$, then $P$ is not in the support of the image of $$\varprojlim R^i\hat{S}(\Omega_e) \to R^i\hat{S}(\Omega_e)$$ for any $e$.
\end{enumerate}

If $\{R^i\hat{S}(\Omega_e)\}$ satisfies the Mittag-Leffler condition for any $i\geq 0$, then (3) also implies (1), (1') and (2) and, moreover, the support of the image of the map in (3) is a closed subset.
\end{thm} \medskip

We also record the following variant of the implication $(2)\Rightarrow (3)$ in the previous theorem.

\begin{prop}
Let $\pi:A \proj W$ be a projection between abelian varieties with generic fiber dimension $f$ and with $\dim W=k$. Let $\{\W_e\}_e$ be a Cartier module on $A$ and let $S_{A,\hat{W}}$ be the Fourier-Mukai functor with kernel $\left(\pi \times 1_{\hat{W}}\right)^{\ast}\PP^{W\times \hat{W}}$. Denote $\L_e = RS_{A,\hat{W}} (D_A (\W_e))$. If $P\in \hat{W}$ is a scheme-theoretic point with $\dim(P)>k+f-\ell$, then $P$ is not in the support of the image of the map $$\varprojlim_e R^{\ell}S_{A,\hat{W}}(\W_e) \lto R^{\ell}S_{A,\hat{W}}(\W_e)$$ Moreover, if the inverse system $\{R^{\ell}S_{A,\hat{W}}(\W_e)\}_e$ satisfies the Mittag-Leffler condition, then the support of the image of the above map is closed and its codimension is $\geq \ell-f$. \label{GV-k}
\end{prop}

\begin{pf}
The proof is identical (modulo shifts) to that of Theorem 4.2 in ~\cite{wz14}, but we include it for the sake of completeness. \medskip

We need to show that if $P\in \hat{W}$ is a scheme-theoretic point with $\dim(P)>k+f-\ell$, then $$\left(R^{\ell}S_{A,\hat{W}}(\varprojlim \W_e)\right)_P \stackrel{0}{\lto} \left(R^{\ell}S_{A,\hat{W}}(\W_e)\right)_P$$

Note in the first place that for any $\ell$, we have the following isomorphisms

\begin{eqnarray} \E xt^p(\L_e,\mathcal{O}_{\hat{W}}) &\simeq&
\mathcal{H}^{p-k}(D_{\hat{W}}(\L_e)) \simeq
\mathcal{H}^{p-k}(D_{\hat{W}}(RS_{A,\hat{W}}(D_A(\W_e)))) \nonumber \\
&\stackrel{[\ast]}{\simeq}& \mathcal{H}^{p}(
R\tilde{S}_{A,\hat{W}}(D_A(D_A(\W_e)))) \simeq
\mathcal{H}^{p}(R\tilde{S}_{A,\hat{W}}(\W_e)) \label{comp1} \end{eqnarray}

and

\begin{eqnarray} \E xt^p(\L,\mathcal{O}_{\hat{W}}) &\simeq&
\mathcal{H}^{p-g}(D_{\hat{W}}(\L)) \simeq
\mathcal{H}^{p-g}(D_{\hat{W}}(\hocolim_e RS_{A,\hat{W}}(D_A(\W_e)))) \nonumber \\
&\simeq& \mathcal{H}^{p-g}(\holim_e D_{\hat{W}}(RS_{A,\hat{W}}(D_A(\W_e)))) \nonumber \\ &\simeq&
\mathcal{H}^p(\holim_e (-1_{\hat{W}})^* RS_{A,\hat{W}}(\tilde{\W}_e)) \nonumber \\
&\simeq&
(-1_{\hat{W}})^*\H^p\left(\holim_e RS_{A,\hat{W}}(\tilde{\W}_e)\right) \label{comp2} \end{eqnarray}

where in $[\ast]$ we used Lemma 2.2 in ~\cite{pp11}. Here, if $S$ is the Fourier-Mukai functor with kernel $\PP$, $\tilde{S}$ denotes the Fourier-Mukai functor with kernel $\PP^{\vee}$; the codimension computation is unaffected by this, so we omit the tildes in the remainder of the proof. \medskip

From the following factorization of the map $R^{\ell}S_{A,\hat{W}}(\varprojlim \W_e) \to R^{\ell}S_{A,\hat{W}}(\W_e)$

$$R^{\ell}S_{A,\hat{W}}(\varprojlim_e \W_e) \lto \H^{\ell}\left( \holim RS_{A,\hat{W}}(\W_e) \right) \lto \H^{\ell}\left( RS_{A,\hat{W}}(\W_e) \right)$$

we see that it suffices to show that the map $$\E xt^{\ell}(\L,\O_{\hat{W}})_P \lto \E xt^{\ell}(\L_e,\O_{\hat{W}})_P$$ is zero for $\dim(P)>k+f-\ell$. In order to see this, we may proceed as in the proof of Theorem 4.2 in ~\cite{wz14}, computing the above map via the commutative diagram $$\xymatrix{ \E xt^i(\H^j(\L),\O_{\hat{W}})_P \ar[d]_{\phi^{i,j}} \ar@{=>}[r] & \E xt^{\ell}(\L,\O_{\hat{W}})_P \ar[d]^{\phi^{\ell}} \\ \E xt^i(\H^j(\L_e),\O_{\hat{W}})_P \ar@{=>}[r] & \E xt^{\ell}(\L_e,\O_{\hat{W}})_P}$$

with $i-j=\ell$. We seek to apply Lemma \ref{spec-seq-zero-map} with $a=\ell-1-f$. If $i>a$, then $i \geq \ell-f > [k+f-\dim(P)] -f = k-\dim(P)$ and hence $\E xt^i(\H^j(\L_e),\O_{\hat{W}})_P=0$ and if $i\leq a$, then $j = i-\ell \leq (\ell -1 -f) -\ell = -1-f$, so that $\H^j(\L)=0$ (c.f. proof of Theorem 3.1.1 in ~\cite{hp13}). Lemma \ref{spec-seq-zero-map} then implies that $\phi^{\ell}=0$ as claimed.
\end{pf} \medskip

\section{Main technical result}

Let $\{\Omega_e\}$ be an inverse system of coherent sheaves on a $g$-dimensional abelian variety satisfying the Mittag-Leffler condition and let $\Omega=\varprojlim \Omega_e$. Let $\Lambda_e=RS_{A,\hat{A}}(D_A(\Omega_e))$ and $\Lambda=\hocolim\Lambda_e$. Suppose that $\{\W_e\}$ is a GV-inverse system, in the sense that $\mathcal{H}^i(\Lambda)=0$ for any $i\neq 0$. It was shown in Theorem 4.2 of ~\cite{wz14} that if $\H^0(\L)$ is torsion-free, then the maps $$\left(\varprojlim R^kS_{A,\hat{A}}(\W_e)\right)_P \lto R^kS_{A,\hat{A}}(\W_e)_P$$ are zero for any point $P\in \hat{A}$ such that $\dim(P) \geq g-k$. We next show a partial converse to this statement. \medskip

In the sequel, we will say that $\H^0(\L)$ has torsion if it is not torsion-free. More concretely, we will say that $\H^0(\L)$ has torsion at a point $P$ if there exists a section $s \in \O_{\hat{A}}$ such that the multiplication map $\H^0(\L)_P \stackrel{\times s_P}{\lto} \H^0(\L)_P$ is not injective. \medskip

Before stating our main result we need to introduce some notation. Consider the following commutative diagram \begin{equation} \xymatrix{0 \ar[r] & \L^t \ar[r] & \L \ar[r] & \FF=\L/\L^t \ar[r] & 0 \\ 0 \ar[r] & \tilde{\L}_e^t \ar[r] \ar[u] & \tilde{\L}_e \ar[r] \ar[u] & \FF_e \ar[u] \ar[r] & 0} \label{main-thm-notation} \end{equation} where $\L^t$ denotes the torsion subsheaf of $\L$, $\tilde{\L}_e=\textrm{Im}(\L_e \to \L)$, $\FF_e=\textrm{Im}(\tilde{\L}_e \to \FF)$ and $\tilde{\L}_e^t=\ker\left( \tilde{\L}_e \to \FF\right)$. It is easy to see that the second row is exact and that $\tilde{\L}_e^t$ is the torsion subsheaf of $\tilde{\L}_e$. It is also clear by construction that $\L = \varinjlim_e \tilde{\L}_e$ (c.f. Theorem \ref{generic-vanishing-char-p}). \medskip

Since the direct limit is exact, by its universal property there is also a commutative diagram $$\xymatrix{0 \ar[r] & \L^t \ar[r] & \L \ar[r] & \FF=\L/\L^t \ar[r] & 0 \\ 0 \ar[r] & \varinjlim_e \tilde{\L}_e^t \ar[r] \ar[u] & \varinjlim_e \tilde{\L}_e \ar[r] \ar[u]_{\simeq} & \varinjlim_e \FF_e \ar[u] \ar[r] & 0}$$ Note that $\L^t \simeq \varinjlim_e \tilde{\L}_e^t$. Indeed, an element $\eta \in \L^t \hookrightarrow \L$ lifts to a class $[\tilde{\eta}_e] \in \varinjlim \tilde{\L}_e = \L$ and this class maps to 0 under the composition $\varinjlim \tilde{\L}_e \to \FF$, so it lies in $\varinjlim \tilde{\L}_e^t$ (by exactness of the direct limit). By the 5-lemma, we have that the right vertical map is also an isomorphism. \medskip

We are now ready to state our main technical result: \medskip

\begin{thm}
Let $\{\Omega_e\}$ be an inverse system of coherent sheaves on a $g$-dimensional abelian variety satisfying the Mittag-Leffler condition. Let $\Lambda_e=RS_{A,\hat{A}}(D_A(\Omega_e))$ and $\Lambda=\hocolim\Lambda_e$. Denote as above $\tilde{\L}_e=\textrm{Im}(\L_e \to \L)$ and define $\tilde{\Omega}_e = RS_{\hat{A},A} (D_{\hat{A}} (\tilde{\L}_e))$. Suppose that $\{\W_e\}$ is a GV-inverse system, in the sense that $\mathcal{H}^i(\Lambda)=0$ for any $i\neq 0$. If $\H^0(\L)$ has a torsion point $P$ of maximal dimension $g-k$, then $$\varprojlim \left(R^kS_{A,\hat{A}}(\tilde{\W}_e)\otimes k(P)\right) \neq 0.$$ Equivalently, the maps $$\varprojlim \left(R^kS_{A,\hat{A}}(\tilde{\W}_e)\otimes k(P)\right) \to R^kS_{A,\hat{A}}(\tilde{\W}_e)\otimes k(P)$$ are non-zero for every $e\gg0$. \label{torsion-non-zero-map}
\end{thm}

\begin{pf}
We start by performing a sequence of reductions. \medskip

\textbf{Reduction 1}: With the notation introduced in diagram (\ref{main-thm-notation}), in order to show that $$\varprojlim \left(R^kS_{A,\hat{A}}(\tilde{\W}_e)\otimes k(P)\right) \neq 0$$ it is sufficient to show that $$\varprojlim_e \left[ \E xt^k(\tilde{\L}_e^t, \O_{\hat{A}}) \otimes k(P) \right] \neq 0.$$

Indeed, consider the long exact sequence for $\E xt$ induced by the short exact sequence $$0 \lto \tilde{\L}_e^t \lto \tilde{\L}_e \lto \FF_e \lto 0$$ namely $$\cdots \lto \E xt^k(\tilde{\L}_e, \O_{\hat{A}}) \lto \E xt^k(\tilde{\L}_e^t, \O_{\hat{A}}) \lto \E xt^{k+1}(\FF_e, \O_{\hat{A}}) \lto \cdots$$ By Lemma 6.3 in ~\cite{pp08b} it follows that $\E xt^{k+1}(\FF_e, \O_{\hat{A}}) \otimes k(P) = 0$ for all $e$, so since $\otimes k(P)$ is right-exact, we obtain a surjection $$\E xt^k(\tilde{\L}_e, \O_{\hat{A}}) \otimes k(P) \lto \E xt^k(\tilde{\L}_e^t, \O_{\hat{A}}) \otimes k(P)$$ and hence\footnote{Note that the system $\left\{ \E xt^k(\tilde{\L}_e^t, \O_{\hat{A}}) \otimes k(P)\right\}_e$ satisfies the ML-condition.} a surjection $$\varprojlim_e \left[ \E xt^k(\tilde{\L}_e, \O_{\hat{A}}) \otimes k(P) \right] \lto \varprojlim_e \left[ \E xt^k(\tilde{\L}_e^t, \O_{\hat{A}}) \otimes k(P) \right].$$

Therefore, if $\varprojlim_e \left[ \E xt^k(\tilde{\L}_e^t, \O_{\hat{A}}) \otimes k(P) \right] \neq0$, then $$\varprojlim_e \left[ \E xt^k(\tilde{\L}_e, \O_{\hat{A}}) \otimes k(P) \right] \stackrel{[\ast]}{\simeq} \varprojlim \left(R^kS_{A,\hat{A}}(\tilde{\W}_e)\otimes k(P)\right) \neq 0$$ as claimed, where the isomorphism $[\ast]$ follows from the following computation

\begin{eqnarray} \E xt^p(\tilde{\L}_e,\mathcal{O}_{\hat{A}}) &\simeq&
\mathcal{H}^{p-g} \left( D_{\hat{A}}(\tilde{\L}_e) \right) \simeq
\mathcal{H}^{p-g} ( (-1_{\hat{A}})^{\ast} RS_{A,\hat{A}} \overbrace{RS_{\hat{A},A} (D_{\hat{A}}(\tilde{\L}_e))}^{:=\tilde{\W}_e} [g]) \nonumber \\
&\simeq& \mathcal{H}^p \left( (-1_{\hat{A}})^* RS_{A,\hat{A}}(\tilde{\W}_e) \right). \end{eqnarray}

\textbf{Reduction 2}: Denoting $i:Z:=\overline{\{P\}}\hookrightarrow \hat{A}$, we next reduce to showing that $$\varprojlim_e \left[ \E xt^k(Li^{\ast}\tilde{\L}^{t}_e,\O_Z) \otimes k(P) \right] \neq0$$

In order to see this, it suffices to show that for every $e$ there is an isomorphism \begin{equation} \E xt^k(\tilde{\L}^t_e,\O_{\hat{A}}) \otimes k(P) \simeq \E xt^k(Li^{\ast}\tilde{\L}_e^t, \O_Z) \otimes k(P). \label{key-iso} \end{equation}

But observe that \begin{multline*} \E xt^k(\tilde{\L}^t_e,\O_{\hat{A}})_{|Z} \otimes k(P) \simeq L^0i^{\ast} \E xt^k(\tilde{\L}^t_e,\O_{\hat{A}}) \otimes k(P) \stackrel{[1]}{\simeq} \H^k\left( Li^{\ast} R\H om(\tilde{\L}^t_e,\O_{\hat{A}}) \right) \otimes k(P) \\ \stackrel{[2]}{\simeq} \E xt^k(Li^{\ast}\tilde{\L}^t_e, \O_Z) \otimes k(P) \end{multline*}

where:

\begin{enumerate}[(i)]
\item The isomorphism in [1] follows from Grothendieck's spectral sequence\footnote{C.f. equation (3.10) in ~\cite{huy06}. Also note that since tensoring by $k(P)$ is exact on $D(Z)$, there is a spectral sequence as written.} $$E_2^{p,q} = L^pi^{\ast}\E xt^q(\tilde{\L}^t_e, \O_{\hat{A}}) \otimes k(P) \Longrightarrow L^{p+q}i^{\ast}R\H om(\tilde{\L}^t_e,\O_{\hat{A}}) \otimes k(P)$$ Note in the first place that $\E xt^q(\tilde{\L}^t_e,\O_{\hat{A}})=0$ near $P$ for all $q<k$. Indeed, since $\tilde{\L}_e^t$ is supported on $Z$\footnote{Note that if $Z$ is an irreducible component of $\Supp \L$, then it is also an irreducible component of $\Supp \H^0(\L_e)$ for every $e>>0$. Let $Z$ be one such component, denote by $\I_Z$ its ideal sheaf and take a section $\eta=\{\eta_e\}\in \L$ supported on $Z$. It is then clear that $Z\subset \Supp \H^0(\L_e)$ for every $e>>0$. Now, if $f\in \I_Z$ is in $\ann(\eta)$, it is clear that $f\in \ann(\eta_e)$ for $e>>0$, so that $Z=\Supp \eta_e$ for all $e>>0$, as claimed.} near $P$, which has codimension $k$, our claim follows from the fact that $\E xt^q(\bullet,\O_{\hat{A}})=0$ for all $q<\codim\Supp(\bullet)$ (c.f. Lemma 6.1 in ~\cite{pp08b}). \medskip

    The differentials coming out of $E_2^{0,k} = L^0i^{\ast}\E xt^k(\tilde{\L}^t_e, \O_{\hat{A}})$ are hence trivial and the differential targeting $E_2^{0,k}$ is $$d_2^{-2,k+1}: L^{-2}i^{\ast}\E xt^{k+1}(\tilde{\L}^t_e, \O_{\hat{A}}) \otimes k(P) \lto L^0i^{\ast}\E xt^k(\tilde{\L}^t_e, \O_{\hat{A}}) \otimes k(P)$$ which is also trivial since there is an open neighborhood $U$ of $P$ over which $$\left[L^{-2}i^{\ast}\E xt^{k+1}(\tilde{\L}^t_e, \O_{\hat{A}})\right]_{|U} \simeq L^{-2} i_{U}^{\ast} \E xt^{k+1}(\tilde{\L}^t_e, \O_{\hat{A}})_{U} = 0$$ where $i_U: Z \cap U \hookrightarrow U$ and where the last vanishing follows from the coherence of $\tilde{\L}^t_e$ and the fact that $P$ has codimension $k$.

\item The isomorphism in [2] follows from the $Lf^{\ast}R\H om(\FF^{\bullet},\mathcal{G}^{\bullet}) \simeq R\H om(Lf^{\ast}\FF^{\bullet}, Lf^{\ast}\mathcal{G}^{\bullet})$ (c.f. equation (3.17) in ~\cite{huy06}).
\end{enumerate}

Finally, in order to show that $$\varprojlim_e \left( \E xt^k(Li^{\ast} \tilde{\L}_e^{t},\O_Z)\otimes k(P) \right)$$ is non-zero, we will use the isomorphism $$\varprojlim_e \left( \E xt^k_{\O_Z} (Li^{\ast} \tilde{\L}_e^{t},\O_Z)\otimes k(P) \right) = \varprojlim_e \E xt^k_{k(P)} \left( Li^{\ast} \tilde{\L}_e^{t} \otimes k(P), k(P) \right) \neq 0$$

where we used that $k(P)\simeq \O_P$ and Proposition III.6.8 in ~\cite{har77}. Denote by $i:Z=\overline{\{P\}} \hookrightarrow \hat{A}$ the inclusion. By Grothendieck duality (c.f. discussion in section 2.3 from ~\cite{bst12}), since all the higher direct images of a closed immersion are zero, we have a functorial isomorphism \begin{equation} R\H om_{\hat{A}} \left( i_{\ast}\left[Li^{\ast}\tilde{\L}_e \otimes k(P) \right], k(P)[-k] \right) \simeq Ri_{\ast} R\H om_{Z} \left( Li^{\ast}\tilde{\L}_e \otimes k(P), Li^! k(P) \right) \label{GD} \end{equation}

Taking k-th cohomology, we obtain \begin{eqnarray}\H om_{k(P)} \left( i_{\ast}\left[ Li^{\ast}\tilde{\L}_e^t \otimes k(P) \right], k(P) \right) &\simeq& \H^k\left(Ri_{\ast}R\H om_{k(P)}(Li^{\ast}\tilde{\L}_e^t \otimes k(P), Li^! k(P))\right) \nonumber \\ &\simeq& i_{\ast}\E xt^k_{k(P)}(Li^{\ast}\tilde{\L}_e^t \otimes k(P), Li^! k(P)). \label{GD2} \end{eqnarray}

Note that the inverse limit of the left hand side is $$\varprojlim \H om_{\hat{A}} \left( i_{\ast}\left[ Li^{\ast}\tilde{\L}_e^t \otimes k(P) \right], k(P) \right) = \H om_{\hat{A}} \left( i_{\ast}\left[ Li^{\ast}\L \otimes k(P) \right], k(P) \right).$$

We claim that the latter sheaf is non-zero. Indeed, note that \begin{equation} H om_{\hat{A}} \left( i^{\ast}\L \otimes k(P) , k(P) \right) \neq0 \label{non-zero-hom} \end{equation}

since it is simply the $k(P)$-dual of the non-zero $k(P)$-vector space $\L_{|Z} \otimes k(P)$. The natural map $Li^{\ast}\L \to L^0 i^{\ast} \L$ induces $$R^0i_{\ast} \left[ Li^{\ast}(\L) \otimes k(P) \right] \stackrel{\simeq}{\lto} R^0i_{\ast} \left[ i^{\ast}(\L) \otimes k(P) \right] \stackrel{\neq0}{\lto} k(P)$$ where the last map is just a non-zero morphism from (\ref{non-zero-hom}) with the source sheaf extended by zero. Since $\L$ is locally free in a neighborhood of $P$ and closed immersions have no higher direct images, the fact that the first map is an isomorphism follows from the degeneration of the spectral sequence (c.f. equation (3.10) in ~\cite{huy06})  $$R^s i_{\ast} \left( L^t i^{\ast} \L \otimes k(P) \right) \simeq R^s i_{\ast} \left( \H^t \left( L i^{\ast} \L \otimes k(P) \right)\right) \stackrel{s+t=p}{\Longrightarrow} R^pi_{\ast} \left( Li^{\ast} \L \otimes k(P) \right)$$ since $L^t i^{\ast} \L \otimes k(P)=0$ for all $t<0$. Hence, the inverse limit on the right hand side of (\ref{GD2}) is also non-zero, and in particular that $$\varprojlim \E xt^k_Z(Li^{\ast}\L_e \otimes k(P), Li^! k(P)) \neq0$$

But recall that $Z$ is a torsion translate of an abelian subvariety of $\hat{A}$, so $\w_Z\simeq \O_Z \simeq \O_P \simeq k(P)$, so $Li^! k(P)\simeq k(P)$ and we may hence conclude that $$\varprojlim \E xt^k_Z(Li^{\ast}\L_e \otimes k(P), k(P)) \neq0$$ as claimed.
\end{pf} \bigskip

\begin{rmk}
Theorem \ref{torsion-non-zero-map} has shown that $$\varprojlim \left(R^k S_{A,\hat{A}}(\tilde{\W}_e)\otimes k(P)\right) \neq 0.$$ where recall, we defined $\tilde{\Omega}_e = RS_{\hat{A},A} (D_{\hat{A}} (\tilde{\L}_e))$ where $\tilde{\L}_e=\textrm{Im}(\L_e \to \L)$. In what follows we will drop the tildes in order to ease de notation: all we need is a projective system of coherent sheaves satisfying the generic vanishing property and inducing a non-zero limit as stated in the theorem.
\end{rmk} \medskip

If $A$ has no super-singular factors, we know by Proposition 3.3.5 in ~\cite{hp13} that for every $e\geq0$, the top dimensional components of the set of points $P\in \hat{A}$ such that the map $\H^0(\L_e)_P\to \H^0(\L)_P$ is non-zero is a torsion translate of an abelian subvariety of $\hat{A}$. \medskip

Let $P\in \hat{A}$ be a torsion point of maximal dimension (namely, $\dim(P)$ is maximal such that $\H^0(\L)_P$ has torsion) and consider $W=\overline{\{P\}}$. In particular, $W$ is a component of $\Supp \H^0(\L)$, and we already argued earlier that $W$ must also be an irreducible component of $\Supp \H^0(\L_e)$ for $e>>0$, so $W$ is also a top dimensional component of the support of the image of the map $\H^0(\L_e)\to \H^0(\L)$. We thus conclude that if $P\in \hat{A}$ is a torsion point of maximal dimension, then $W=\overline{\{P\}}$ is a torsion translate of an abelian subvariety of $\hat{A}$. \medskip

In this context, Theorem \ref{torsion-non-zero-map} yields the following.

\begin{corol}
Let $\{\Omega_e\}$ be a Mittag-Leffler inverse system of coherent sheaves on a
$g$-dimensional abelian variety $A$ with no supersingular factors, and let $\Omega=\varprojlim \Omega_e$. Let
$\Lambda_e=RS_{A,\hat{A}}(D_A(\Omega_e))$ and $\Lambda=\hocolim\Lambda_e$.
Suppose that $\{\W_e\}$ is a GV-inverse system, in the sense that $\mathcal{H}^i(\Lambda)=0$ for any $i\neq 0$. If $P$ is a torsion point of $\H^0(\L)$ of maximal dimension (so that $W=\overline{\{P\}}$ is a torsion translate of an abelian subvariety of $\hat{A}$), then there are non-zero maps $$\varprojlim \left(R^{g-k}S_{A,\hat{W}}(\W_e)\otimes k(P)\right) \to R^{g-k}S_{A,\hat{W}}(\W_e)\otimes k(P)$$ for  $e\gg0$, where the Fourier-Mukai kernel of $S_{A,\hat{W}}$ is given by $\PP^{A\times \hat{W}}=(id\times \iota)^{\ast}\PP^{A\times \hat{A}}$, $\PP^{A\times \hat{A}}$ being the normalized Poincar{\'e} bundle of $A\times \hat{A}$. \label{torsion-non-zero-map2}
\end{corol}

\begin{pf}
By Theorem \ref{torsion-non-zero-map} we have a non-zero map $$\varprojlim \left(R^{g-k}S_{A,\hat{A}}(\W_e)\otimes k(P)\right) \to R^{g-k}S_{A,\hat{A}}(\W_e)\otimes k(P)$$ for some $e>0$. Consider the base-change maps $$R^{g-k}S_{A,\hat{A}}(\W_e)\otimes k(P) \to H^{g-k}(A,\W_e\otimes \PP^{A\times \hat{A}}_{|A\times \{P\}}), \quad R^{g-k}S_{A,\hat{W}}(\W_e)\otimes k(P) \to H^{g-k}(A,\W_e\otimes \PP^{A\times \hat{W}}_{|A\times \{P\}})$$ Note that the second map is an isomorphism by flat base change: indeed, denoting by $\iota:\hat{W} \hookrightarrow \hat{A}$ the inclusion, we have \begin{eqnarray*} R^{g-k}S_{A,\hat{W}}(\W_e)\otimes k(P) &\stackrel{def}{=}& R^{g-k}p_{\hat{W}\ast}\left(p_A^{\ast}\W_e \otimes \PP^{A\times \hat{W}} \right)\otimes k(P) \\ &\stackrel{FBC}{\simeq}& R^{g-k}p_{\hat{W}\ast}'\left(p_A^{\ast}\W_e \otimes \PP^{A\times \hat{W}}\otimes k(P) \right) \\ &\stackrel{def}{\simeq}& R^{g-k}p_{\hat{W}\ast}'\left(p_A^{\ast}\W_e \otimes (id\times \iota)^{\ast}\PP^{A\times \hat{A}}\otimes k(P) \right) \\ &\stackrel{[\ast]}{\simeq}& H^{g-k}(A,\W_e \otimes \PP^{A\times \hat{W}}_{|A\times \{P\}}) \end{eqnarray*} where in $[\ast]$ we used Proposition III.8.5 in ~\cite{har77} and where $p_{\hat{W}}'$ is the base change of the projection, as illustrated in the diagram $$\xymatrix{A \times_{\hat{W}} \Spec k(P) \ar[r] \ar[d]_{p_{\hat{W}}'} & A \times \hat{W} \ar[d]^{p_{\hat{W}}} \\ \Spec k(P) \ar[r]^{\textrm{flat}} & \hat{W}}$$

However note that we have \begin{equation} H^{g-k}(A,\W_e\otimes \PP^{A\times \hat{A}}_{|A\times \{P\}}) \simeq H^{g-k}(A,\W_e\otimes \PP^{A\times \hat{W}}_{|A\times \{P\}}) \label{obvious-iso} \end{equation} so we can write both base change maps in the following diagram $$\xymatrix{\varprojlim\left(R^{g-k}S_{A,\hat{A}}(\W_e)\otimes k(P)\right) \ar[r]^-{\neq 0} \ar[d] & R^{g-k}S_{A,\hat{A}}(\W_e)\otimes k(P) \ar[d]^-{[\ast]} \\ \varprojlim H^{g-k}\left(A,\W_e\otimes \PP^{A\times \hat{A}}_{|A\times \{P\}}\right) \ar[r] \ar[d]_{\simeq} & H^{g-k}\left(A,\W_e\otimes \PP^{A\times \hat{A}}_{|A\times \{P\}}\right) \ar[d]^{\simeq} \\ \varprojlim H^{g-k}\left(A,\W_e\otimes \PP^{A\times \hat{W}}_{|A\times \{P\}}\right) \ar[r] & H^{g-k}\left(A,\W_e\otimes \PP^{A\times \hat{W}}_{|A\times \{P\}}\right) \\ \varprojlim\left(R^{g-k}S_{A,\hat{W}}(\W_e)\otimes k(P)\right) \ar[u]^{\simeq} \ar[r] & R^{g-k}S_{A,\hat{W}}(\W_e)\otimes k(P) \ar[u]_{\simeq}}$$ where the top horizontal map is non-zero for $e\gg0$ by Theorem \ref{torsion-non-zero-map}, the middle isomorphisms are the ones in (\ref{obvious-iso}), and where the isomorphisms at the bottom follow from flat base change as described above. \medskip

We seek to show that the bottom horizontal map is non-zero for some $e$. Nevertheless, note that if it this were not the case, then the top horizontal map could not possibly be non-zero, since in any case the base change maps $[\ast]$ are injective by the proof of Proposition III.12.5 in ~\cite{har77}\footnote{In a nutshell, let $f:X\to Y=\Spec A$ be a projective morphism and let $\FF$ be a coherent sheaf on $X$. For any $A$-module $M$, define $T^i(M):=H^i(X,\FF\otimes_A M)$, which is a covariant additive functor from $A$-modules to $A$-modules which is exact in the middle (by Proposition III.12.1 in ~\cite{har77}). Writing $$A^r\to A^s \to M \to 0$$ we have a diagram $$\xymatrix{T^i(A)\otimes A^r \ar[r] \ar[d]_{\simeq} & T^i(A)\otimes A^s \ar[r] \ar[d]_{\simeq} & T^i(A) \otimes M \ar[r] \ar[d]^{\varphi} & 0 \\ T^i(A^r) \ar[r] & T^i(A^s) \ar[r] & T^i(M)}$$ where $\varphi:T^i(A)\otimes M \to T^i(M)$ is the base change map and where the two first vertical arrows are isomorphisms. A straight-forward diagram chase then shows that $\varphi$ is injective.}.
\end{pf} \medskip

\begin{rmk}
We will be using two different Fourier-Mukai kernels on $A\times \hat{W}$. If $\iota:\hat{W}\hookrightarrow \hat{A}$ denotes the inclusion and $\pi:A\proj W$ is the dual projection, we have a diagram $$\xymatrix{A\times \hat{W} \ar[r]^{id_A\times \iota} \ar[d]_{\pi \times id_{\hat{W}}} & A \times \hat{A} \\ W \times \hat{W} & }$$ If $\PP^{A\times \hat{A}}$ and $\PP^{W\times \hat{W}}$ denote the normalized Poincar{\'e} bundles on $A\times \hat{A}$ and $W\times \hat{W}$ respectively, on $A\times \hat{W}$ we may consider the locally-free sheaves $(id_A\times \iota)^{\ast}\PP^{A\times \hat{A}}$ and $(\pi\times id_{\hat{W}})^{\ast}\PP^{W\times \hat{W}}$. \medskip

In Corollary \ref{torsion-non-zero-map2} we proved a non-vanishing statement for the derived Fourier-Mukai transform with respect to the former kernel and in what follows we need an analogous statement for the transform with respect to the latter. Nevertheless, note that we are simply looking at fibers over points $w\in \hat{W}\subset \hat{A}$ (concretely over the generic point of $\hat{W}$), and over these points both sheaves are isomorphic. Indeed, $w\in \hat{W}\subset \hat{A}$ determines $\PP^{W\times \hat{W}}_{|W\times \{w\}}\in \Pic(W)$ and $\PP^{A\times \hat{A}}_{|A\times \{w\}}\in \Pic(A)$, with $\PP^{A\times \hat{A}}_{|A\times \{w\}}\simeq \pi^{\ast} \PP^{W\times \hat{W}}_{|W\times \{w\}}$, and therefore $$\overbrace{\left[(id_A\times \iota)^{\ast}\PP^{A\times \hat{A}}\right]_{|A\times \{w\}}}^{\simeq \PP^{A\times \hat{A}}_{|A\times \{w\}}} \simeq \overbrace{\left[(\pi\times id_{\hat{W}})^{\ast}\PP^{W\times \hat{W}}\right]_{|A\times \{w\}}}^{\simeq \pi^{\ast} \PP^{W\times \hat{W}}_{|W\times \{w\}}}$$ \label{different-FM-kernels}
\end{rmk} \medskip

\section{Fibering of the Albanese image}

In ~\cite{el97}, Ein and Lazarsfeld showed the following statement:

\begin{thm}[see ~\cite{el97}, Theorem 3]
$X$ is a smooth projective variety of maximal Albanese dimension over a field of characteristic zero and $\chi(\w_X)=0$, then the image of the Albanese map is fibered by subtori of $A$. \label{EL-fibered-by-tori}
\end{thm}

\begin{skpf} We sketch the proof given in ~\cite{pp08} (Theorem E). \medskip
If $\chi(\w_X)=0$, it follows that $V^0(\w_X)\subset \hat{A}$ is a proper subvariety (c.f. Lemma 1.12(b) in ~\cite{par11}). Choose an irreducible component $W\subset V^0(\w_X)$ of codimension $p>0$, which is a torsion translate of an abelian subvariety of $\hat{A}$ that we also denote by $W$. Let $\pi:A\proj \hat{W}$ denote the dual projection and consider the diagram $$\xymatrix{X \ar[r]^-a & Y:=a(X) \ar[d]_-{h=\pi_{|Y}} \ar@{^{(}->}[r] & A \ar@{->>}[dl]^{\pi} \\ & \hat{W}}$$

Since the fibers of the projection $A\proj \hat{W}$ are abelian subvarieties of dimension $p$, the conclusion of the theorem will follow provided that $f\geq p$, where $f$ denotes the dimension of the generic fiber of $h$. In order to see this, recall the following standard facts:

\begin{enumerate}[(a)]
\item $a_{\ast}\w_X$ is a GV-sheaf on $Y=a(X)$ and $V^0(\w_X)=V^0(a_{\ast}\w_X)$.

\item If $W$ is an irreducible component of $V^0(a_{\ast}\w_X)$ of codimension $p$, then it is also a component of $V^p(a_{\ast}\w_X)$.

\item $a_{\ast}\w_X$ is a $GV_{-f}$-sheaf with respect the the Fourier-Mukai functor with kernel $(\pi \times 1_W)^{\ast}\PP^{W\times \hat{W}}$, so in particular $\codim_W V^p(a_{\ast}\w_X) \geq p-f$ for every $p\geq0$.
\end{enumerate} \medskip

By (b) we have $W \subseteq V^p(a_{\ast}\w_X) \subseteq W$ so that $\codim_W V^p(a_{\ast}\w_X)=0$, and finally (c) yields $f\geq p$, which is what we sought to show.
\end{skpf} \medskip

Our goal in this section is to prove a positive characteristic analogue of Theorem \ref{EL-fibered-by-tori}. \medskip

Let $X$ be a smooth projective variety of maximal Albanese dimension and denote by $a:X\to A$ the Albanese map. Then $a_{\ast}\w_X$ is a Cartier module and we may consider the inverse system $\{\W_e=F_{\ast}^e S^0a_{\ast}\w_X\}_e$. Define $\L_e=RS_{A,\hat{A}}D_A(\W_e)$ and set $\L=\hocolim \L_e$. By Corollary 3.3.1 in ~\cite{hp13} we have that $H^i(\W\otimes P_{\alpha})=0$ for every $i>0$ and very general $\alpha\in \hat{A}$. Thus, defining as above $\chi(\W):=\chi(\W\otimes P_{\alpha})$ for very general $\alpha\in \hat{A}$, we see that $\chi(\W)=h^0(\W\otimes P_{\alpha})$ and it seems that in trying to extend Theorem \ref{EL-fibered-by-tori} to positive characteristic, one should assume that $h^0(\W\otimes P_{\alpha})=0$. \medskip

This leaves us in a setting which is similar to the one we encountered in the proof of Theorem \ref{Main-Theorem}. If $rk(\L)=h^0(\W\otimes P_{\alpha})=0$, then in particular $\Lambda$ must be a torsion sheaf. In light of this observation, we show the following:

\begin{thm}
Let $X$ be a smooth projective variety of maximal Albanese dimension and let $a:X\to A$ be a generically finite map to an abelian variety $A$ with no supersingular factors. Let $g=\dim A$. Consider the inverse system $\{\W_e=F_{\ast}^eS^0a_{\ast}\w_X\}_e$ and denote $\W=\varprojlim \W_e$. Define $\L_e=RS_{A,\hat{A}}D_A(\W_e)$ and assume that the sheaf $\H^0(\L)=\varinjlim \H^0(\L_e)$  has torsion. Then the image of the Albanese map is fibered by linear subvarieties of $\hat{A}$. \label{alb-fibered-by-tori}
\end{thm}

\begin{pf}
Let $w\in \hat{A}$ be a torsion point of $\H^0(\L)$ of maximal dimension $k$; by the remark preceding Corollary \ref{torsion-non-zero-map2}, we have that $\hat{W}:=\overline{\{w\}}\subset \hat{A}$ is a torsion translate of an abelian subvariety of $\hat{A}$, which we still denote by $\hat{W}$. Denote by $\pi:A\to W$ the projection dual to the inclusion $\hat{W}\hookrightarrow \hat{A}$. $$\begin{diagram} \node{Y:=a(X)} \arrow{e,t,J}{} \arrow{s,l}{h=\pi_{|Y}} \node{A} \arrow{sw,r}{\pi} \\ \node{W} \end{diagram}$$

By Corollary \ref{torsion-non-zero-map2} we know that the map $$\varprojlim \left(R^{g-k}S_{A,\hat{W}}(\W_e)\otimes k(w)\right) \to R^{g-k}S_{A,\hat{W}}(\W_e)\otimes k(w)$$ is non-zero for every $e\gg 0$, where the Fourier-Mukai kernel of $S_{A,\hat{W}}$ is given by $(id\times \iota)^{\ast}\PP$, $\PP$ being the normalized Poincar{\'e} bundle of $A\times \hat{A}$. \medskip

Recall that, in general, even though the system $\{\W_e\}$ satisfies the Mittag-Leffler condition, the inverse system $\{R^tS(\W_e)\}_e$ may fail to do so (c.f. Example 3.2 in ~\cite{wz14}). We handle the Mittag-Leffler case first, however, since the proof is neater and the subsequent generalization does not rely on new ideas. \medskip

\textbf{Case in which $\{R^{g-k}S_{A,\hat{W}}(\W_e)\}_e$ satisfies the ML-condition}. The proof in this case goes along the lines of that of Theorem E in ~\cite{pp08}. Note that if $w$ is the generic point of $\hat{W}\hookrightarrow \hat{A}$, we have the following:

\begin{eqnarray*} w &\stackrel{[1]}{\in}& \left\{w\in \hat{W}: \quad \varprojlim_e \left(R^{g-k}S_{A,\hat{W}}(\W_e) \otimes k(w)\right) \stackrel{\neq0}{\lto} R^{g-k}S_{A,\hat{W}}(\W_e) \otimes k(w) \right\} \\ &\stackrel{[2]}{=}& \left\{ w \in \hat{W}: \quad \left(\varprojlim_e R^{g-k}S_{A,\hat{W}}(\W_e)\right) \otimes k(w) \stackrel{\neq0}{\lto} R^{g-k}S_{A,\hat{W}}(\W_e) \otimes k(w) \right\} \\ &\subseteq& \left\{ w \in \hat{W}: \quad \left(\varprojlim_e R^{g-k}S_{A,\hat{W}}(\W_e)\right)_w \stackrel{\neq0}{\lto} \left(R^{g-k}S_{A,\hat{W}}(\W_e)\right)_w \right\} \\ &\subseteq& \hat{W} \end{eqnarray*}

where $w$ lies in the first set by Corollary \ref{torsion-non-zero-map2} and where the equality [2] follows from Lemma \ref{inverse-limit-tensor-product-commute}, since we are under the assumption that the system $\{R^{g-k}S_{A,\hat{W}}(\W_e)\}_e$ satisfies the Mittag-Leffler condition. \medskip

This implies that the codimension (in $\hat{W}$) of the support of image of the map $$\varprojlim_e R^{g-k}S_{A,\hat{W}}(\W_e) \lto R^{g-k}S_{A,\hat{W}}(\W_e)$$ is zero (this support is closed - under the Mittag-Leffler assumption - by Proposition 4.3 in ~\cite{wz14}). At the same time, by Proposition \ref{GV-k} we know that this codimension must be $\geq g-k-f$, where $f$ is the dimension of a general fiber of $h$, so in particular $f\geq g-k$ and this concludes the proof under the Mittag-Leffler assumption (indeed, the fibers of the projection $A\proj W$ are abelian subvarieties of dimension $g-k$). \medskip

\textbf{General case}. We finally observe that, in our setting, we can actually do without the Mittag-Leffler assumption on $\{R^{g-k}S_{A,\hat{W}}(\W_e)\}_e$. We only used this assumption in order to guarantee the closedness of the support of the image of the map $\left(\varprojlim_e R^{g-k}S_{A,\hat{W}}(\W_e)\right)_w \stackrel{\neq0}{\lto} \left(R^{g-k}S_{A,\hat{W}}(\W_e)\right)_w$ and in order to ensure that the inverse limit commutes with $\otimes k(w)$. \medskip

Note in the first place that we don't need the support of the image of the above map to be closed for the previous argument to work. Proposition \ref{GV-k} shows that in order for $w$ to belong to the support, we need $\codim \overline{\{w\}}\geq g-k-f$, and this suffices in order to conclude that $f\geq g-k$. \medskip

With regards to the commutation of the inverse limit and $\otimes k(w)$, note in the first place that there is always an inclusion $\supseteq$ induced by the natural map $$\varprojlim_e \left( R^{g-k}S_{A,\hat{W}}(\W_e) \right) \otimes k(w) \lto \varprojlim_e \left(R^{g-k}S_{A,\hat{W}}(\W_e) \otimes k(w)\right)$$ Moreover, in our setting, the opposite inclusion $\subseteq$ follows from the flatness of $k(w)$ as an $\O_{\hat{W}}$-module and the fact that the projection formula and its consequences still hold in the category of quasi-coherent sheaves under some perfection assumptions (c.f. Lemma 71 in ~\cite{mur06}). We state this below as a lemma, and the proof is hence complete.
\end{pf} \medskip

\begin{lemma}
With the same notations as above, if $$\varprojlim_e \left(R^{g-k}S_{A,\hat{W}}(\W_e) \otimes k(w)\right) \stackrel{\neq0}{\lto} R^{g-k}S_{A,\hat{W}}(\W_e) \otimes k(w)$$ then $$\varprojlim_e \left( R^{g-k}S_{A,\hat{W}}(\W_e) \right) \otimes k(w) \stackrel{\neq0}{\lto} R^{g-k}S_{A,\hat{W}}(\W_e) \otimes k(w)$$ In particular, equality [2] in the above chain still holds.
\end{lemma}

\begin{pf}
As in the proof of Proposition \ref{GV-k}, let $\tilde{\L}_e = R\tilde{S}_{A,\hat{W}}(D_A(\W_e))$, where $\tilde{S}_{A,\hat{W}}$ denotes the Fourier-Mukai transform with kernel $\mathcal{P}^{\vee}$, with $\mathcal{P}=\left(\pi \times 1_{\hat{W}}\right)^{\ast}\PP^{W\times \hat{W}}$. Note in the first place that we have the following isomorphisms of $\O_{\hat{W}}$-modules: \begin{eqnarray*} \varprojlim_e \left(R^{g-k}S_{A,\hat{W}}(\W_e) \otimes k(w)\right) &\stackrel{[1]}{\simeq}& \varprojlim_e \left(\E xt^{g-k}_{\O_{\hat{W}}}(\tilde{\L}_e,\O_{\hat{W}}) \otimes k(w)\right) \\ &\simeq& \varprojlim_e \left(\H^{g-k}\left(R\H om_{\O_{\hat{W}}} (\tilde{\L}_e,\O_{\hat{W}}) \right) \otimes k(w)\right) \\ &\stackrel{[2]}{\simeq}& \varprojlim_e \H^{g-k}\left(R\H om_{\O_{\hat{W}}} (\tilde{\L}_e,\O_{\hat{W}}) \otimes k(w) \right) \\ &\stackrel{[3]}{\simeq}& \H^{g-k}\left(\varprojlim_e \left(R\H om_{\O_{\hat{W}}} (\tilde{\L}_e,\O_{\hat{W}}) \otimes k(w) \right)\right) \\ &\stackrel{[4]}{\simeq}& \H^{g-k}\left(\varprojlim_e R\H om_{\O_{\hat{W}}} (\tilde{\L}_e,k(w)) \right) \\ &\stackrel{[5]}{\simeq}& \H^{g-k}\left(R\H om_{\O_{\hat{W}}} (\tilde{\L},k(w)) \right) \\ &\stackrel{[6]}{\simeq}& \H^{g-k}\left(R\H om_{\O_{\hat{W}}} (\tilde{\L},\O_{\hat{W}}) \otimes k(w) \right) \\ &\stackrel{[7]}{\simeq}& \H^{g-k}\left(R\H om_{\O_{\hat{W}}} (\tilde{\L},\O_{\hat{W}})\right) \otimes k(w) \\ &\stackrel{[8]}{\simeq}& \H^{g-k}\left(\holim RS_{A,\hat{W}}(\W_e)\right) \otimes k(w) \end{eqnarray*}

where [1] and [8] follow from the computations (\ref{comp1}) and (\ref{comp2}) in the proof of Proposition \ref{GV-k}, [2] and [7] follow from the flatness of $\otimes k(w)$ as an $\O_{\hat{W}}$-module, [3] follows from Proposition \ref{inverse-limits-commute-functor}, since the system $\left\{\E xt^{p-1}(\L_e,\O_{\hat{W}})\otimes k(w)\right\}_e$ satisfies the ML-condition, and [4] and [6] follow from the isomorphism\footnote{This isomorphism holds for complexes of sheaves of modules $\mathcal{F},\mathcal{G},\mathcal{H}$ provided that either $\mathcal{F}$ or $\mathcal{H}$ are perfect (c.f. Lemma 71 in ~\cite{mur06}). Note that $k(w)$ is a perfect complex, being a coherent sheaf.} $$R\H om(\FF,\mathcal{G}) \otimes \mathcal{H} \simeq R\H om(\FF,\mathcal{G} \otimes \mathcal{H}).$$

Hence, by assumption we have a non-zero map $$\H^{g-k}\left(\holim RS_{A,\hat{W}}(\W_e)\right) \otimes k(w) \stackrel{\neq0}{\lto} R^{g-k}S_{A,\hat{W}}(\W_e) \otimes k(w)$$ and the conclusion of the lemma then follows from the following commutative diagram $$\xymatrix{\varprojlim_e \left(R^{g-k}S_{A,\hat{W}}(\W_e) \right) \otimes k(w) \ar[r] & R^{g-k}S_{A,\hat{W}}(\W_e) \otimes k(w) \\ \H^{g-k}\left(\holim RS_{A,\hat{W}}(\W_e)\right) \otimes k(w) \ar[u] \ar[ur]_{\neq0}}$$
\end{pf} \medskip

In particular, within the context of principally polarized abelian varieties, the same argument yields the following statement:

\begin{corol}
Let $(A,\Theta)$ is a principally polarized abelian variety with no supersingular factors defined over an algebraically closed field of characteristic $p>0$. Assume further that $\Theta$ is irreducible. Consider the inverse system $\{\W_e:=F_{\ast}^e(\w_{\Theta}\otimes \tau_{\Theta})\}_e$ on $A$ and set $\L=\hocolim_e RS_{A,\hat{A}}D_A(\W_e)$. Then $\L$ is a torsion-free quasi-coherent sheaf concentrated in degree 0. \label{theta-divisors-not-ruled}
\end{corol}

\begin{pf}
The fact that $\L=\H^0(\L)$ is a quasi-coherent sheaf concentrated in degree zero follows from Theorem \ref{generic-vanishing-char-p}(i), since $\w_{\Theta}$ is a Cartier module. \medskip

Assume for a contradiction that $\H^0(\L)$ is not torsion-free and fix an irreducible component $\hat{W}:=\overline{\{w\}}\hookrightarrow \hat{A}$ of maximal dimension of the closure of the set of torsion points of $\H^0(\L)$. Denote by $\pi:A\proj W$ the dual projection. $$\xymatrix{\Theta \ar@{^{(}->}[r] \ar[d]_{h} & A \ar@{->>}[dl]_{\pi} \\ W & }$$

We may then argue as in the proof of Theorem \ref{alb-fibered-by-tori} to conclude that $\Theta$ is fibered by abelian subvarieties of $A$, but this is not possible given that $\Theta$ is irreducible (and hence of general type) in light of Abramovich's work (c.f. Section 2.3 or ~\cite{abr95}).
\end{pf} \medskip

\section{Singularities of Theta divisors}

We now focus on the singularities of Theta divisors and embark on the proof of Theorem \ref{Main-theorem}. As a warm-up, we focus on simple abelian varieties to start with, namely those which do not contain smaller dimensional abelian varieties.

\subsection{Case of simple abelian varieties}

The crux of the argument resides in the construction of sections of $\O_A(\Theta)$ which vanish along the test ideal $\tau(\Theta)$ and, in the case of simple abelian varieties, it is a direct consequence of the results in ~\cite{hp13}.  The proof of the general case will follow the same pattern, albeit further work will be required to prove that the required sections exist. \medskip

\begin{thm}
Let $(A,\Theta)$ be a PPAV over an algebraically closed field $K$ of characteristic $p>0$ such that $A$ is simple and ordinary. Then $\Theta$ is strongly F-regular. (In particular, $\Theta$ is F-rational, and by Lemma 2.34 in ~\cite{bst12}, it is normal and Cohen-Macaulay).
\label{Main-theorem-simple}
\end{thm} \medskip

\begin{pf}
On $\hat{A}$ consider the inverse system $\W_e=F_{\ast}^e\W_0$, where $\W_0=\omega_{\Theta}\otimes \tau(\Theta)$. This yields a direct system $\Lambda_e= R\hat{S}D_A\W_e$ equipped with natural maps $\H^0(\Lambda_e) \to \H^0(\Lambda)=\Lambda=\hocolim \Lambda_e$. By \ref{generic-vanishing-char-p}, we know that $\Lambda$ is quasi-coherent sheaf in degree 0. \medskip

Consider the set $$Z=\left\{\alpha\in \hat{A}: \quad Im \left(\H^0(\Lambda_0)\otimes \O_{\hat{A},\alpha} \lto \Lambda \otimes \O_{\hat{A},\alpha} \right) \neq 0\right\}$$ By Proposition \ref{GV-corollary2} we know that $Z$ is a finite union of torsion translates of subtori. Since $\hat{A}$ is simple by assumption, this implies that either $Z=\hat{A}$ or $Z$ is a finite set. \medskip

Assume for a contradiction that $Z$ is finite. By Proposition \ref{GV-corollary3} we have $t_x^{\ast}\Omega=\Omega$ for every $x\in \widehat{\hat{A}/Z}=\widehat{\hat{A}}$, so that $\Supp \Omega=A$. Since the maps in the inverse system $F_{\ast}^e\Omega_0=F_{\ast}^e(\omega_{\Theta}\otimes \tau(\Theta))$ are surjective, we know by Proposition \ref{GV-corollary3} that $\Supp \Omega=\Supp \Omega_0=\Supp \w_{\Theta}\otimes \tau(\Theta)=A$, which is absurd. \medskip

We must thus have $Z=\hat{A}$, so that $\Supp \H^0(\Lambda_0)=\hat{A}$, and hence cohomology and base change yields $H^0(A,\w_{\Theta}\otimes \tau(\Theta) \otimes P_{\alpha})\neq 0$ for all $\alpha\in \hat{A}$. \medskip

Consider the following commutative diagram: $$\xymatrix{ H^0(A,P_{\alpha}) \ar[r] & H^0(A,\O_A(\Theta)\otimes P_{\alpha}) \ar[r] & H^0(\Theta,\O_A(\Theta)_{|\Theta} \otimes P_{\alpha}) \\ H^0(A,K\otimes P_{\alpha}) \ar[r] \ar@{^{(}->}[u] & H^0(A,\O_A(\Theta)\otimes P_{\alpha} \otimes \tilde{\tau}) \ar[r] \ar@{^{(}->}[u] & H^0(\Theta,\O_A(\Theta)_{|\Theta}\otimes \tau(\Theta) \otimes P_{\alpha}) \ar@{^{(}->}[u] & }$$

where $K$ is defined so that the diagram commutes. In the top row we have $H^1(A,P_{\alpha})=0$ for $\alpha\neq 0$ since $P_{\alpha}$ is topologically trivial. The polarization induced by $\Theta$ is principal, so $h^0(\O_A(\Theta)\otimes P_{\alpha})=1$. Since by the above discussion $H^0(\Theta,\O_A(\Theta)_{|\Theta}\otimes P_{\alpha}\otimes \tau(\Theta))\neq 0$, it follows that $H^0(\Theta,\O_A(\Theta)_{|\Theta}\otimes P_{\alpha})\neq 0$ and hence both the right inclusion and the top right restriction are equalities. \medskip

The ideal sheaf $\tilde{\tau}$ on $A$ is defined as follows: fix an open subset $U=\Spec R\subseteq A$ and assume that $\Theta$ is given by an ideal sheaf $I=I(\Theta)$. Let $J=\tau_{\Theta}(U)$ be the test ideal of $\Theta$ and let $\tilde{J}\subset R$ be an ideal such that $J=\tilde{J}/I$. Omitting the twist by $P_{\alpha}$, the diagram above locally boils down to

$$\xymatrix{ && R/\tilde{J} \ar[r]^-{\simeq} & (R/I)/(\tilde{J}/I) \simeq R/\tilde{J} & \\ 0 \ar[r] & I \ar[r] & R \ar[r] \ar[u] & R/I \ar[r] \ar[u] & 0 \\ 0 \ar[r] & I \ar[r] \ar[u] & \tilde{J} \ar[r] \ar[u] & J \ar[r] \ar[u] & 0 }$$

so $\tilde{\tau}(U)=\tilde{J}$. Now taking cohomology, a section $s\in H^0(J)$ embeds as $\bar{s}\in H^0(R/I)$ and maps to zero in $H^0(R/\tilde{J})$ by exactness. By exactness of the second row, $\bar{s}$ lifts to $\tilde{\bar{s}}\in H^0(R)$, which still projects to zero in $H^0(R/\tilde{J})$ by commutativity of the top square, so $\tilde{\bar{s}}$ must lift to a non-zero section of $H^0(\tilde{J})$. \medskip

Finally, since $H^0(\Theta,\O_A(\Theta)_{|\Theta}\otimes P_{\alpha}\otimes \tau(\Theta))\neq 0$ for every $\alpha\in \hat{A}$ and these sections lift to sections of $H^0(A,\O_A(\Theta)\otimes P_{\alpha})$ vanishing along $\tilde{\tau}$, we conclude that $h^0(A,\O_A(\Theta)\otimes P_{\alpha}\otimes \tilde{\tau})=1$. Hence if $\tilde{\tau}$ were not trivial, we would have $Zeros(\tilde{\tau})\subset \Theta+\alpha_P$ for every $\alpha_P\in A$ (where $\alpha_P\in A$ is the point corresponding to $P_{\alpha}\in \Pic^0(A)$), which is absurd since these translates of $\Theta$ don't have any points in common. We thus conclude that $\tilde{\tau}=\O_A$, and hence $\tau(\Theta)=\O_{\Theta}$ so that $\Theta$ is strongly F-regular.
\end{pf} \medskip

In the proof of Theorem \ref{Main-theorem-simple} we used the simplicity of $A$ in order to show that $H^0(A,\w_{\Theta}\otimes \tau(\Theta) \otimes P_{\alpha})\neq 0$ for all $\alpha\in \hat{A}$. The same argument we employed above will work in the general case provided that we can show the existence of sections in $H^0(A,\w_{\Theta}\otimes \tau(\Theta) \otimes P_{\alpha})$, and it turns out that this is somewhat more involved.

\subsection{General case}

We finally study singularities of Theta divisors in the general setting. As we mentioned earlier, the argument will be analogous to the one employed to prove the theorem in the case of simple abelian varieties, although additional work will be required to prove that there exist sections in $H^0(\Theta,\O_A(\Theta)_{|\Theta} \otimes \tau(\Theta) \otimes P_{\alpha})$. \medskip

The main ingredient in Ein and Lazarsfeld's proof over fields of characteristic zero was that given a smooth projective variety $X$ of maximal Albanese dimension such that $\chi(X,\w_X)=0$, the image of its Albanese morphism is fibered by tori (c.f. Theorem 3 in ~\cite{el97}). Our proof will rely on Corollary \ref{theta-divisors-not-ruled}, where we proved that if the sheaf $\H^0(\L)$ associated to the inverse system $\{F_{\ast}^e(\w_{\Theta}\otimes \tau_{\Theta})\}$ was not torsion-free, then $\Theta$ would be fibered by tori, which is impossible since $\Theta$ is irreducible. \medskip

In a nutshell, and as in the case of simple abelian varieties, $\Theta$ will be strongly F-regular provided that there exist non-trivial sections in $H^0(\Theta,\O_A(\Theta)\otimes \O_{\Theta} \otimes \tau(\Theta) \otimes P_{\alpha})$ and we will show that if that was not the case, then the sheaf $\H^0(\L)$ would have torsion, a contradiction. \medskip

\begin{thm}
Let $(A,\Theta)$ be an ordinary principally polarized abelian variety over an algebraically closed field $k$ of characteristic $p>0$. If $\Theta$ is irreducible, then $\Theta$ is strongly F-regular. \label{Main-Theorem}
\end{thm}

\begin{pf}
The proof goes along the lines of Theorem \ref{Main-theorem-simple}: consider again the commutative diagram: $$\xymatrix{ H^0(A,P_{\alpha}) \ar[r] & H^0(A,\O_A(\Theta)\otimes P_{\alpha}) \ar[r] & H^0(\Theta,\O_A(\Theta)_{|\Theta} \otimes P_{\alpha}) \\ H^0(A,K\otimes P_{\alpha}) \ar[r] \ar@{^{(}->}[u] & H^0(A,\O_A(\Theta)\otimes P_{\alpha} \otimes \tilde{\tau}) \ar[r] \ar@{^{(}->}[u] & H^0(\Theta,\O_A(\Theta)_{|\Theta}\otimes \tau(\Theta) \otimes P_{\alpha}) \ar@{^{(}->}[u] & }$$

In the proof of Theorem \ref{Main-theorem-simple} we used the simplicity of $A$ to conclude easily that $$H^0(\Theta,\overbrace{\O_A(\Theta)\otimes \O_{\Theta}}^{\w_{\Theta}} \otimes \tau(\Theta) \otimes P_{\alpha})\neq 0$$ It then followed from the commutative diagram above that $H^0(A,\O_A(\Theta)\otimes P_{\alpha} \otimes \tilde{\tau})\neq 0$ and this in turn forced $\tilde{\tau}$ to be trivial, whence $\tau(\Theta)=\O_{\Theta}$. \medskip

We shall now use the previous results in order to conclude that $H^0(\Theta,\w_{\Theta} \otimes \tau(\Theta) \otimes P_{\alpha})\neq 0$. In a nutshell, assuming for a contradiction that $\tau(\Theta)$ is not trivial we will show that $0\neq S^0(\Theta,\w_{\Theta} \otimes \tau(\Theta) \otimes P_{\alpha}) \subseteq H^0(\Theta,\w_{\Theta} \otimes \tau(\Theta) \otimes P_{\alpha})$ for general $\alpha\in \hat{A}$. The diagram above then yields $H^0(A,\O_A(\Theta)\otimes P_{\alpha} \otimes \tilde{\tau})\neq 0$ and since by assumption $\tau(\Theta)\ne \O_{\Theta}$, we conclude that $\textrm{Zeroes}(\tilde{\tau})\subset \Theta+\alpha_P$ for general $\alpha_P\in A$ (as before, $\alpha_P$ denotes the point in $A$ corresponding to $P_{\alpha}\in \Pic^0A$), but this is not possible since general translates of $\Theta$ do not have points in common. Therefore we must have $\tau(\Theta)=\O_{\Theta}$, and hence $\Theta$ is strongly F-regular. \medskip

We now argue by contradiction: let $\W=\varprojlim F_{\ast}^eS^0\w_{\Theta}= \varprojlim F_{\ast}^e(\w_{\Theta}\otimes \tau(\Theta))$. By Corollary 3.2.1 in ~\cite{hp13} and its proof, for all closed $\alpha\in \hat{A}$ we have that $$\H^0(\L)\otimes k(\alpha) \simeq H^0(\Theta,\varprojlim \W_e \otimes P_{\alpha}^{\vee})^{\vee}$$ so assuming for a contradiction that $S^0(\Theta,\w_{\Theta} \otimes \tau(\Theta) \otimes P_{\alpha})=0$ for general $\alpha\in \hat{A}$, it follows that $rk (\Lambda)=0$, and hence that $\H^0(\Lambda)$ has torsion. Nevertheless, this is not possible by Corollary \ref{theta-divisors-not-ruled}, since we are assuming that $\Theta$ is irreducible, and this concludes the proof.
\end{pf}


\small
\begin{thebibliography}{XXX99}
\bibitem[Abr95]{abr95} D. Abramovich, \emph{Subvarieties of semiabelian varieties}, 1995.
\bibitem[BL92]{bl92} C. Birkenhake and H. Lange, \emph{Complex abelian varieties}, Vol. 302. Springer, 1992.
\bibitem[BS12]{bs12} M. Blickle and K. Schwede, \emph{$p^{-1}$-linear maps in algebra and geometry}, arXiv:1205.4577v2, 2012
\bibitem[BST12]{bst12} M. Blickle, K. Schwede and K. Tucker, ~\emph{F-singularities via alterations}, arXiv:1107.3807, 2012
\bibitem[Car08]{car08} J. Carter, \emph{The Morava K-Theory Eilenberg-Moore spectral sequence}, New York J. Math. \textbf{14}, 495-515, 2008
\bibitem[EL97]{el97} L. Ein and R. Lazarsfeld \emph{Singularities of theta divisors and the birational geometry of irregular varieties}, J. Amer. Math. Soc \textbf{10}:1, 243-258, 1997
\bibitem[GL87]{gl87} M. Green and R. Lazarsfeld, \emph{Deformation theory, generic vanishing theorems and some conjectures of Enriques, Catanese and Beauville}, Invent. Math. \textbf{90}, 389-407, 1987
\bibitem[GL91]{gl91} M. Green and R. Lazarsfeld, \emph{Higher obstructions to deforming cohomology groups of line bundles}, Jour. of. A.M.S. \textbf{4}, 87-103, 1991
\bibitem[Hac04]{hac04} C. Hacon. \emph{A derived category approach to generic vanishing}, J. Reine Angew. Math. \textbf{575}, 173-187, 2004
\bibitem[Hac11]{hac11} C. Hacon. \emph{Singularities of pluri-theta divisors in characteristic $p>0$}, arXiv:1112.2219v1, 2011
\bibitem[HH90]{hh90} M. Hochster and C. Huneke. \emph{Tight closure, invariant theory, and the Briançon-Skoda theorem}, J. Amer. Math. Soc. 3 (1990), no. 1, 31-116
\bibitem[HK12]{hk12} C. Hacon and S. Kov{\'a}cs. \emph{Generic vanishing fails for singular varieties and in characteristic $p>0$}, arXiv:1212.5105
\bibitem[HP13]{hp13} C. Hacon and Z. Patakfalvi. \emph{Generic vanishing in characteristic $p>0$ and the characterization of ordinary abelian varieties}, arXiv:1310.2996, 2013
\bibitem[Har77]{har77} R. Hartshorne. \emph{Algebraic geometry.} Graduate Texts in Mathemaics, \textbf{52}, Springer, 1977.
\bibitem[Har78]{har78} R. Hartshorne. \emph{On the de Rham cohomology of algebraic varieties}, Publ. Math. de l'IHES, vol. \textbf{45}, 5-99, 1978, 2007
\bibitem[Har98]{har98} N. Hara. \emph{A characterization of F-rational singularities in terms of injectivity of Frobenius maps}, American Journal of Mathematics, Volume 120, Number 5, pp. 981-996, 1998
\bibitem[HW02]{hw02} N. Hara and K.I. Watanabe. \emph{F-regular and F-pure rings vs. log terminal and log canonical singularities}, J. Algevraic Geom. \textbf{11}, no. 2, 363-392, 2002
\bibitem[Huy06]{huy06} D. Huybrechts, \emph{Fourier-Mukai transforms in algebraic geometry}, Oxford Mathematical Monographs, 2006
\bibitem[KV81]{kv81} Y. Kawamata and E. Viehweg. \emph{On a characterization of abelian varieties in the classification theory of algebraic varieties}, Comp. Math. \textbf{41}, 355-360, 1981
\bibitem[Laz04]{laz04} R. Lazarsfeld. \emph{Positivity in Algebraic Geometry I and II}, \textbf{48-49}, Springer-Verlag, 2004
\bibitem[Mor00]{mor00} S. Mori. \emph{Classification of higher-dimensional algebraicd varieties.} Proc. Symp. Pue Math. \textbf{46} (1987), 269-332
\bibitem[Muk81]{muk81} S. Mukai. \emph{Duality between $D(X)$ and $D(\hat{X})$ with its application to Picard scheaves}, Nagoya Math. J. \textbf{81}, 133-175, 1981
\bibitem[Mur06]{mur06} D. Murfet. \emph{Derived categories of quasi-coherent sheaves}, 2006
\bibitem[Nee96]{nee96} A. Neeman. \emph{The Grothendieck duality theorem via Bousfield's techniques and Brown representability}, J. Amer. Math. Soc. \textbf{9}, 205-235, 1996
\bibitem[OSS80]{oss80} C. Okonek, M. Schneider and H. Spindler. \emph{Vector bundles on complex projective spaces}, Progress in Mathematics \textbf{3}, Birkhäiser, Boston, 1980.
\bibitem[Par11]{par11} G. Pareschi. \emph{Basic results on irregular varieties via Fourier-Mukai methods}, Current Developments in Algebraic Geometry, MRSI Publications, \textbf{59}, 2011
\bibitem[PP03]{pp03} G. Pareschi and M. Popa, \emph{Regularity on abelian varieties I}, arXiv:math/0110003
\bibitem[PP08]{pp08} G. Pareschi and M. Popa, \emph{Regularity on abelian varieties III: relationship with generic vanishing and apploications}, arXiv:0802.1021, 2008
\bibitem[PP08b]{pp08b} G. Pareschi and M. Popa, \emph{Strong generic vanishing and a ahigher dimensional Castelnuovo-de Franchis inequality}, arXiv:0808.2444, 2008
\bibitem[PP11]{pp11} G. Pareschi and M. Popa. \emph{GV-sheaves, Fourier-Mukai transform and generic vanishing}, Smer. J. Math. J. \textbf{133}:1, 235-271, 2011
\bibitem[PR03]{pr03} R. Pink and D. Roesddler. \emph{A conjecture of Beauville and Catanese revisited}, Mathematische Annalen, Volume 330, Issue 2, pp 293-308, 2004
\bibitem[Sch09]{sch09} K. Schwede. \emph{F-adjunction}, arXiv:0901.1154
\bibitem[Sch12]{sch12} K. Schwede. \emph{A canonical linear system associated to adjoint divisors in characteristic $p>0$}, arXiv:1107.3833
\bibitem[Smi97]{smi97} K. Smith. \emph{F-rational rings have rational singularities}, Amer. J. Math. 119 (1997), no. 1, 159–180, 1997.
\bibitem[Tak04]{tak04} S. Takagi. \emph{An interpretation of multiplier ideals via tight closure}, J. Algebraic. Geom. \textbf{13}, 393–415, 2004
\bibitem[Wan11]{wan11} J. Wang. \emph{Quotients of algebraic groups}, 2011
\bibitem[Wei94]{wei94} C. A. Weibel. \emph{An introduction to homological algebra}, Cambridge University Press, 1994
\bibitem[WZ14]{wz14} A. Watson and Y. Zhang. \emph{On the generic vanishing theorem of Cartier modules}, arXiv:1404.2669
\end{thebibliography}
\end{document}